\documentclass[a4paper,11pt,column]{quantumarticle}
\pdfoutput=1

\usepackage[utf8]{inputenc}
\usepackage[english]{babel}
\usepackage[T1]{fontenc}

\usepackage{amsmath,amssymb,amsfonts,bm,dsfont}
\usepackage{physics}
\usepackage{breqn}
\usepackage{booktabs}
\usepackage[numbers,sort&compress]{natbib}
\usepackage{hyperref}
\usepackage{tikz}
\usepackage{tikz-cd}
\usetikzlibrary{patterns,arrows}
\usepackage{float}
\usepackage{algpseudocode}
\usepackage{multirow}
\usepackage{xcolor}

\makeatletter
\newcounter{algorithm}
\renewcommand{\thealgorithm}{\arabic{algorithm}}
\newcommand{\fps@algorithm}{tbp}
\newcommand{\ftype@algorithm}{4}
\newcommand{\ext@algorithm}{loa}
\newcommand{\fnum@algorithm}{Algorithm~\thealgorithm}

\makeatother
\renewcommand{\vec}[1]{\bm{#1}}
\begin{document}

\title{A tensor-train multidimensional inverse Laplace transform}

\author{Martin Mikkelsen}
\email{Martin.Mikkelsen@di.ku.dk}
\affiliation{Department of Computer Science, University of Copenhagen}
\author{Michael Kastoryano}
\affiliation{Department of Computer Science, University of Copenhagen}

\maketitle

\begin{abstract}
Laplace transforms and their numerical inverses arise throughout applied mathematics, physics, finance, and probability theory. Numerical inversion, however, quickly becomes intractable in high dimensions because the number of quadrature evaluations grows exponentially with dimension. We develop a tensor train (TT) formulation of the multidimensional inverse Laplace transform. The method constructs a TT approximation of the transformed function on the complex quadrature grid and then performs the inversion through a sequence of tensor contractions. Under suitable low-rank assumptions, this reduces the computational cost from exponential to polynomial in the dimension, provided that the relevant bond dimensions remain bounded. The method has only a small number of tunable parameters and admits error estimations. We demonstrate its performance in numerical experiments, benchmarked against Monte Carlo estimates and exact references, for multivariate normal-inverse Gaussian, Wishart, and correlated Gamma-type distributions.

\end{abstract}

\section{Introduction}

Integral transforms are fundamental tools in applied mathematics and computational science. Among them, the Laplace transform is especially important in probability theory, physics, mathematical finance, and the analysis of differential equations. Many applications require solving the inverse problem: recovering a function, or a probability density, from its Laplace transform. Analytic inversion is possible only in special cases, so practical applications often rely on numerical inversion.

Even in one dimension, numerical Laplace inversion is non-trivial: quadrature-based methods require careful contour placement and summation acceleration to control aliasing and truncation errors \cite{hurwitz1956numerical, MultidimensionalLaplace}. Existing approaches include the Post--Widder formula \cite{EpostGeneralized, DVWidder}, series-expansion methods \cite{chung1986taylor, Weeksmethod}, quadrature methods \cite{salzer1955orthogonal, schmittroth1960numerical, duffy1993numerical}, rational-approximation methods \cite{longman1972numerical, sidi1976best}, and Talbot-type methods \cite{talbot1979accurate, green1955calculation, evans1993numerical, de1995parallel}. These methods are well established for one-dimensional transforms, and several have been extended to multidimensional settings \cite{moorthy1995inversion,abate1998numerical,singhal2005numerical}. For a comprehensive review, see \cite{cohen2007numerical}.

 Many applications in finance, probability theory, and stochastic modelling involve multivariate Laplace transforms whose inversion requires high-dimensional oscillatory complex sums or integrals. Direct numerical methods typically scale exponentially with dimension, making them computationally infeasible beyond two or three dimensions.

Tensor-network methods offer a promising route around this bottleneck. Originally developed in condensed-matter physics \cite{perez2006matrix} and independently in numerical analysis \cite{tt_decomp}, they provide efficient low-rank representations of high-dimensional tensors and operators \cite{hubig2017generic}. Among these representations, the tensor-train decomposition---also known as matrix product states in the physics literature---has proved especially useful for high-dimensional numerical problems. TT methods have been applied successfully in uncertainty quantification \cite{zhang2014enabling}, parametric PDEs \cite{richter2021solving, dolgov2015polynomial, HamiltonJacobiBellman, arenstein2025fast, arenstein2025full}, and machine learning \cite{wahls2014learning, sozykin2022ttopt, batsheva2023protes, novikov2015tensorizing}.

Tensor-train and quantics tensor-train methods have also been applied to integral transforms \cite{QFT1, QFT2, Dolgov2012, jaseem2026quantum}, primarily for Fourier transforms or forward Laplace transforms. Here we instead address the inverse problem, which requires evaluating the transformed function on complex contour nodes and poses a qualitatively different multidimensional challenge.

The central observation of this work is that multidimensional Laplace inversion has a natural tensor-network formulation. The inversion formula of Choudhury, Lucantoni, and Whitt \cite{MultidimensionalLaplace} can be expressed as a recursive composition of one-dimensional inversion operators acting independently along each coordinate. We approximate the transformed function $\tilde{f}$ on a quadrature grid by a low-rank tensor train and then contract the Euler quadrature weights along each dimension. When the transformed function has low rank on this complex quadrature grid, the resulting algorithm reduces both storage and computational complexity from exponential to polynomial scaling in the dimension.

The main contributions are threefold. First, we derive a tensor-train representation of the multidimensional inverse Laplace transform, giving an operator formulation suited to low-rank approximation. Second, we propose a numerical inversion algorithm that combines this operator structure with established tensor-train techniques, in particular \texttt{TT-cross} interpolation \cite{oseledets_TT_cross,Submatrix,FastAdaptive,maxvol}. Third, we demonstrate the method on multivariate normal-inverse Gaussian, Wishart, and correlated Gamma-type models, and we study how the accuracy depends on the tensor-train rank and on the inversion parameters.

The remainder of the paper is organized as follows. Section~\ref{sec:prelim} reviews the multidimensional Laplace transform and the inversion method of Choudhury, Lucantoni, and Whitt \cite{MultidimensionalLaplace}. Section~\ref{sec:tt} reviews tensor trains and tensor cross interpolation. Section~\ref{sec:tensor_formulation} presents the tensor-train formulation of the inverse Laplace transform. Section~\ref{sec:error_bounds} discusses error sources and parameter selection. Section~\ref{sec:examples} presents numerical examples and applications, and Section~\ref{sec:distributional_queries} shows how the recovered TT density can be used for marginal probabilities, and conditional distributional queries.

\section{Preliminaries}\label{sec:prelim}

\subsection{The Laplace transform and its inverse}

We consider the unilateral Laplace transform of a function $f(t)$ defined by
\begin{equation}\label{eq:Laplace1D}
    \mathcal{L}\{f(t)\}(s)
    =
    \int_0^\infty e^{-st} f(t)\,dt,
\end{equation}
and denote the transform by $\tilde f(s)$. We assume that $f(t)$ is defined on $[0,\infty)$ and that the integral converges for $\Re(s)>0$.

The multidimensional Laplace transform of a function $f(\vec t)$ with $\vec t=(t_1,\dots,t_d)$ is defined as
\begin{equation}\label{eq:LaplaceMultiDim}
    \mathcal{L}\{f(\vec t)\}(\vec s)=\int_{[0,\infty)^d} e^{-\vec s\cdot \vec t} f(\vec t)\,d\vec t,
\end{equation}
where $\vec s\cdot \vec t$ denotes the Euclidean inner product. As in the one-dimensional case, suitable growth and regularity assumptions are required to ensure convergence and uniqueness.

We denote the inverse Laplace transform by
\begin{equation}\label{eq:LaplaceInvMultiDim}
    f(\vec t)=\mathcal{L}^{-1}\{\tilde f(\vec s)\}.
\end{equation}

In practice, numerical inversion of the Laplace transform is challenging even in one dimension. In multiple dimensions the problem is further complicated by the exponential growth in the number of required function evaluations. This makes direct multidimensional inversion prohibitively expensive in most settings.

We follow the multidimensional inversion framework of \cite{MultidimensionalLaplace}, which treats both continuous variables, via Laplace transforms, and discrete variables, via generating functions. This paper focuses on continuous variables, although the tensor-network formulation extends naturally to the discrete case.
The inversion algorithm of \cite{MultidimensionalLaplace} proceeds by successively applying one-dimensional inversion operators along each coordinate, from $k=d$ down to $k=1$. At each step, an Euler-summation inversion is applied in dimension $k$, while the remaining coordinates are held fixed through their associated quadrature nodes. This yields a recursive operator formulation in which the multidimensional inverse is obtained by composition of one-dimensional operators.

For each dimension $k$, define the quadrature node
\begin{equation}\label{eq:quad_nodes}
    \xi_k(p_k,j_k;t_k)=\frac{A_k}{2t_k\ell_k}-\frac{i p_k\pi}{t_k\ell_k}-\frac{i j_k\pi}{t_k},
\end{equation}
where $p_k=1,\ldots,\ell_k$ and $j_k\in\mathbb Z$. We write
\begin{equation}
    \bm{\xi}(\bm p,\bm j;\bm t)=\bigl(\xi_1(p_1,j_1;t_1),\ldots,\xi_d(p_d,j_d;t_d)\bigr).
\end{equation}
The one-dimensional inversion operator in coordinate $k$ is
\begin{equation}\label{eq:invmultidim_op}
\begin{aligned}
   (\Lambda_k g)(\bm{\xi})
   &:=
   \frac{e^{A_k/(2\ell_k)}}{2t_k\ell_k}
   \sum_{p_k=1}^{\ell_k}
   \sum_{j_k=-\infty}^{\infty}
   (-1)^{j_k} e^{-i p_k\pi/\ell_k}  \\
   &\quad \times
   g\!\left(
      \xi_1,\ldots,
      \xi_k(p_k,j_k;t_k),
      \ldots,\xi_d
   \right),
\end{aligned}
\end{equation}
where the other coordinates are kept fixed during the $k$-th operation. The multidimensional approximation is then obtained by composition,
\begin{equation}\label{eq:invmultidim}
    \bar f(\bm t)=(\Lambda_1\circ\Lambda_2\circ\cdots\circ\Lambda_d)\tilde f,
\end{equation}
where we apply from $k=d$ to $k=1$ and denote the recovered function by $\bar{f}$. 

After all $d$ steps, the output $\bar f(\bm t)$ approximates $f(\bm t)$, with a discretization error $\bar e$ discussed in Section~\ref{sec:error_bounds}. The infinite sums over $j_k$ are truncated to
$j_k=-J,\ldots,J$, with $J=m_E+n_E$, and Euler-accelerated. This gives the finite quadrature grid used below. A full example for $d=2$ is provided in Appendix~\ref{app:2D}. Table~\ref{tab:inversion_parameters} summarizes the inversion parameters used throughout.

\begin{table}[t]
\centering
\small
\begin{tabular}{p{0.25\columnwidth}p{0.65\columnwidth}}
\hline
Parameter & Meaning \\
\hline
$A_k$ & Damping parameter which controls the discretization error \\
$\ell_k$ & Period-splitting parameter which controls the round-off error \\
$n_{\mathrm{E}}$ & Euler truncation parameter \\
$m_{\mathrm{E}}$ & Euler acceleration order \\
$J=n_{\mathrm{E}}+m_{\mathrm{E}}$ & Maximum $j$-index \\
$N_j=2J+1$ & Number of $j$-grid nodes \\
$n_t$ & Number of time-grid points per dimension \\
\hline
\end{tabular}
\caption{Inversion parameters used in the tensor train inversion algorithm.}
\label{tab:inversion_parameters}
\end{table}

\section{Tensor trains}\label{sec:tt}

Tensor trains (TT) provide compressed representations of high-dimensional tensors
\cite{tt_decomp, schollwock2011density}. Instead of storing all
$\prod_{k=1}^d n_k$ entries of a $d$-dimensional tensor explicitly, the TT format
represents the tensor as a sequence of low-dimensional cores whose contraction
reconstructs the full tensor. Specifically, an order-$d$ tensor
$a \in \mathbb{C}^{n_1 \times \cdots \times n_d}$ is expressed as a chain of
third-order TT cores $A_k \in \mathbb{C}^{\alpha_{k-1} \times n_k \times \alpha_k}$,
using the \{left, physical, right\} index convention.

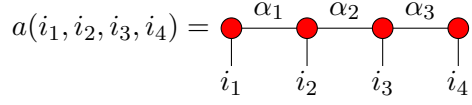
\begin{figure}[ht]
  \centering
  \begin{tikzpicture}
    \node[anchor=east] at (-0.15,0.0) {$a(i_1,i_2,i_3,i_4)=$};

    \draw[black] (0,0) -- node [label={[shift={(0,-0.15)}]$\alpha_1$}] {} ++ (1,0);
    \draw[black] (1,0) -- node [label={[shift={(0,-0.15)}]$\alpha_2$}] {} ++ (1,0);
    \draw[black] (2,0) -- node [label={[shift={(0,-0.15)}]$\alpha_3$}] {} ++ (1,0);

    \draw[black] (0,0) -- node [label={[shift={(0,-0.9)}]$i_1$}] {} ++ (0,-0.5);
    \draw[black] (1,0) -- node [label={[shift={(0,-0.9)}]$i_2$}] {} ++ (0,-0.5);
    \draw[black] (2,0) -- node [label={[shift={(0,-0.9)}]$i_3$}] {} ++ (0,-0.5);
    \draw[black] (3,0) -- node [label={[shift={(0,-0.9)}]$i_4$}] {} ++ (0,-0.5);

    \node[draw,shape=circle,fill=red,scale=0.75] at (0,0){};
    \node[draw,shape=circle,fill=red,scale=0.75] at (1,0){};
    \node[draw,shape=circle,fill=red,scale=0.75] at (2,0){};
    \node[draw,shape=circle,fill=red,scale=0.75] at (3,0){};
  \end{tikzpicture}
  \caption{Tensor network diagram of a tensor train with $d=4$. Each node
    represents a TT core $A_k$, internal edges represent bond dimensions
    $\alpha_k$, and external edges correspond to physical indices $i_k$.
    Boundary bond dimensions $\alpha_0 = \alpha_4 = 1$ are omitted for clarity.}
  \label{fig:mps_ttv}
\end{figure}

In the TT format, the entries of $a$ are written as
\begin{align}\label{eq:ttv}
    a(i_1,\dots,i_d)
    &= \prod_{k=1}^d A_k(:,i_k,:) \\
    &= \sum_{\alpha_1,\ldots,\alpha_{d-1}}
       \prod_{k=1}^d A_k(\alpha_{k-1}, i_k, \alpha_k),
\end{align}
where $A_k(:,i_k,:) \in \mathbb{C}^{\alpha_{k-1}\times\alpha_k}$ is the matrix
slice of core $A_k$ at fixed physical index $i_k$, the integers $\alpha_k$ are
the bond dimensions with boundary values $\alpha_0 = \alpha_d = 1$, and the sum
runs over all bond indices. A tensor
network diagram is shown in Figure~\ref{fig:mps_ttv}.

Each node in Figure~\ref{fig:mps_ttv} represents a TT core $A_k$ with physical index $i_k$ (downward edge) and bond dimensions $\alpha_{k-1}$, $\alpha_k$ (horizontal edges). The maximum bond dimension $\chi = \max_k \alpha_k$ is a measure of the degree of compression, as the storage cost of the TT representation scales as
\begin{equation}\label{eq:storage}
    \sum_{k=1}^{d} \alpha_{k-1}\, n_k\, \alpha_k \;\leq\; d\chi^2 n, \quad n = \max_k n_k,
\end{equation}
compared to $\prod_{k=1}^d n_k$ for the full tensor. Since the right-hand side of equation \eqref{eq:storage} is linear in $d$ for fixed $\chi$ and $n$, the TT format avoids the curse of dimensionality whenever $\chi$ remains bounded as $d$ increases.

\subsection{Tensor cross interpolation}\label{sec:TCI}
An efficient way to construct a tensor-train representation of a multivariate function is the \texttt{TT-cross} algorithm \cite{oseledets_TT_cross}. Consider a function $f: \mathbb{R}^d \rightarrow \mathbb{R}$ discretized on $n$ grid points per dimension. Storing all entries of the resulting tensor requires $n^d$ values, which becomes infeasible even for moderate $d$. Given pointwise access to $f$, however, one can construct a tensor-train approximation using $\mathcal{O}(d n \chi^2)$ function evaluations \cite{vysotsky2021tensor}, where $\chi$ is the bond dimension of the resulting TT.

The key building block is the matrix cross approximation. A matrix $M$ of size $n \times m$ and rank $r$ satisfies
\begin{equation}
    M \approx M(:,J)\,M(I,J)^{-1}\,M(I,:),
\end{equation}
where $I \subset \{1,\ldots,n\}$ and $J \subset \{1,\ldots,m\}$ are index sets of $r$ rows and $r$ columns respectively, $M(:,J)$ is the $n \times r$ submatrix of selected columns, $M(I,:)$ is the $r \times m$ submatrix of selected rows, and $M(I,J)$ is the $r \times r$ intersection submatrix. When $r$ equals the rank of $M$ and the pivots are chosen optimally, the approximation is exact; the computational benefit is largest when $r \ll m, n$.

When the rank of $M$ exceeds $r$, a quasi-optimal cross approximation is obtained by choosing pivot sets $(I,J)$ that maximize $|\det M(I,J)|$, the volume of the intersection submatrix. We find such pivots using the \texttt{MaxVol} algorithm \cite{Submatrix, goreinov2001maximal, maxvol}: starting from a random column set $J$, we select the $r$ rows $I \subset \{1,\ldots,n\}$ of $M(:,J)$ that maximize $|\det M(I,J)|$; using these rows, we then select columns $J$ that maximize the same volume. This alternating procedure is repeated until convergence.

This approach generalizes to the order-$d$ tensor case as follows. At each bond $k$, the tensor $a$ is unfolded into a matrix along the boundary between modes $1,\ldots,k$ and modes $k+1,\ldots,d$. A set of left multi-indices $I_k$ and right multi-indices $J_k$ plays the role of the row and column pivot sets from the matrix case, enumerating configurations of the left and right groups of modes respectively. Sweeping through the TT cores one at a time, the matrix cross approximation is applied to each unfolding while all other cores are held fixed, and the resulting factors update the adjacent cores. The sweeps are repeated until the relative error on a random validation set falls below the tolerance $\varepsilon$, or until a bond dimension reaches the prescribed maximum $\chi_{\text{max}}$.

\section{Tensor network formulation of multidimensional inversion}\label{sec:tensor_formulation}

The inversion formula~\eqref{eq:invmultidim} has a key structural property: each quadrature node $\xi_k$ depends only on the local indices $(p_k,j_k,t_k)$ of dimension $k$. The dimensions are coupled only through the transformed function $\tilde{f}$, evaluated at the combined node $(\xi_1,\ldots,\xi_d)$. This separability makes a tensor-train representation of $\tilde{f}$ well suited to the inversion, which can be split into two stages.
A naive evaluation of \eqref{eq:invmultidim} requires $\mathcal{O}((\ell J)^d)$ evaluations of $\tilde{f}$, where $\ell$ is the number of quadrature points per dimension and $J$ is the range of the Euler index. For example, with $d=5$ and $\ell J=250$, direct summation requires roughly $10^{12}$ evaluations. If $\tilde{f}$ instead admits a low-rank tensor-train approximation on the joint index space $(p_1,j_1,t_1,\ldots,p_d,j_d,t_d)$, the full sum can be evaluated by sequential tensor contractions at cost $\mathcal{O}(d\ell J\chi^2)$, linear in $d$ for fixed bond dimension $\chi$.

The algorithm consists of two parts:
\begin{enumerate}
    \item Build a $3d$-dimensional TT approximation of $\tilde{f}$ on the extended
          quadrature grid $(p_k,j_k,t_k)_{k=1}^d$ using TT-cross interpolation. Each physical dimension contributes three tensor cores. The resulting TT representation is shown in Figure~\ref{fig:tt-structure}.
    \item Contract the $(p_k,j_k)$ legs of each dimension with the Euler weights
          and absorb the time-domain prefactor into the $t_k$ core, following equation \eqref{eq:invmultidim_op}. This produces a $d$-dimensional tensor train whose physical indices are the time variables $(t_1,\ldots,t_d)$, as shown in Figure~\ref{fig:contraction_scheme}.
\end{enumerate}

\subsection{The extended quadrature grid and the \texorpdfstring{$3d$}{3d}-dimensional TT}

\begin{figure}[ht]
  \centering
  \begin{tikzpicture}
    \node[anchor=east] at (-0.15,0.0) {$\tilde{f}=$};

    \draw[dashed,gray!60,rounded corners=2pt] (-0.2,0.28) rectangle (1.6,-1.2);
    \node[above,font=\small] at (0.7,0.28) {$d=1$};

    \draw[dashed,gray!60,rounded corners=2pt] (1.9,0.28) rectangle (3.7,-1.2);
    \node[above,font=\small] at (2.8,0.28) {$d=2$};

    \draw (0,0) -- (0.7,0);
    \draw (0.7,0) -- (1.4,0);
    \draw (1.4,0) -- (2.1,0);
    \draw (2.1,0) -- (2.8,0);
    \draw (2.8,0) -- (3.5,0);
    \draw (3.5,0) -- (3.85,0);
    \draw[dotted] (3.85,0) -- (4.55,0);
    \draw (4.55,0) -- (4.9,0);
    \draw (4.9,0) -- (5.6,0);
    \draw (5.6,0) -- (6.3,0);

    \draw (0,0)   -- ++(0,-0.65) node[below,font=\small]{$p_1$};
    \draw (0.7,0)  -- ++(0,-0.65) node[below,font=\small]{$j_1$};
    \draw (1.4,0)  -- ++(0,-0.65) node[below,font=\small]{$t_1$};
    \draw (2.1,0)  -- ++(0,-0.65) node[below,font=\small]{$p_2$};
    \draw (2.8,0)  -- ++(0,-0.65) node[below,font=\small]{$j_2$};
    \draw (3.5,0)  -- ++(0,-0.65) node[below,font=\small]{$t_2$};
    \draw (4.9,0)  -- ++(0,-0.65) node[below,font=\small]{$p_d$};
    \draw (5.6,0)  -- ++(0,-0.65) node[below,font=\small]{$j_d$};
    \draw (6.3,0)  -- ++(0,-0.65) node[below,font=\small]{$t_d$};

    \foreach \x in {0, 0.7, 1.4, 2.1, 2.8, 3.5, 4.9, 5.6, 6.3}
      \node[draw,circle,fill=red,scale=0.75] at (\x,0){};

  \end{tikzpicture}
  \caption{Structure of the $3d$-dimensional TT for $\tilde{f}$. Cores are grouped in triplets, one per dimension $k$, with physical indices $(p_k, j_k, t_k)$ (dashed boxes). Intra-triplet bonds capture correlations among the quadrature indices within a single dimension; inter-triplet bonds (between boxes) carry correlations across dimensions. Boundary bonds $r_0 = r_{3d} = 1$ are omitted.}
  \label{fig:tt-structure}
\end{figure}
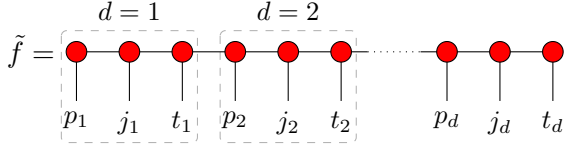

We treat the quadrature indices of each dimension as independent physical modes of a $3d$-dimensional tensor. We use shared parameters $A_k = A$ and $\ell_k = \ell$ across all dimensions. The quadrature nodes in equation~\eqref{eq:quad_nodes} for dimension $k$ then take the form
\begin{equation}
    \xi_k(p_k, j_k, t_k)
    =
    \frac{A}{2t_k\ell}
    -
    \frac{ip_k\pi}{t_k\ell}
    -
    \frac{ij_k\pi}{t_k},
\end{equation}
where $p_k \in \{1,\ldots,\ell\}$, $j_k \in \{-J,\ldots,J\}$, and $t_k$ ranges over a prescribed evaluation grid of $n_t$ points. The $3d$-dimensional tensor
\begin{equation}\label{eq:calF}
    \mathcal{F}(p_1,j_1,t_1,\ldots,p_d,j_d,t_d)
    =
    \tilde{f}\!\bigl(\xi_1,\,\ldots,\,\xi_d\bigr)
\end{equation}
records all evaluations of $\tilde{f}$ on the quadrature grid. We approximate $\mathcal{F}$ in TT format using TT-cross interpolation~\cite{oseledets_TT_cross,maxvol}. The TT built from cross interpolation is arranged in $d$ triplets of tensor cores, one per dimension $k$, with physical legs $(p_k, j_k, t_k)$, as shown in Figure~\ref{fig:tt-structure}. The bond dimensions $r_{3k-2}$ and $r_{3k-1}$ capture correlations among the quadrature indices within a single dimension, while bond dimensions $r_{3k}$ carry correlations between different dimensions.

\subsection{Contracting the quadrature indices}

Using the TT representation of $\tilde{f}(\xi_1,\ldots,\xi_d)$, we evaluate the inversion sum~\eqref{eq:invmultidim} by contracting the quadrature legs $p_k$ and $j_k$ against Euler weight vectors for each dimension $k$.

We name the three TT cores of triplet $k$ as $P^{(k)}$, $J^{(k)}$, and $T^{(k)}$, with shapes $r_{3k-3} \times \ell \times r_{3k-2}$, $r_{3k-2} \times (2J+1) \times r_{3k-1}$, and $r_{3k-1} \times n_t \times r_{3k}$, respectively. Fixing the physical index of a core to a particular value yields a matrix: for example, $P^{(k)}[:,p,:] \in \mathbb{C}^{r_{3k-3} \times r_{3k-2}}$ is the matrix slice of $P^{(k)}$ at quadrature point $p$.

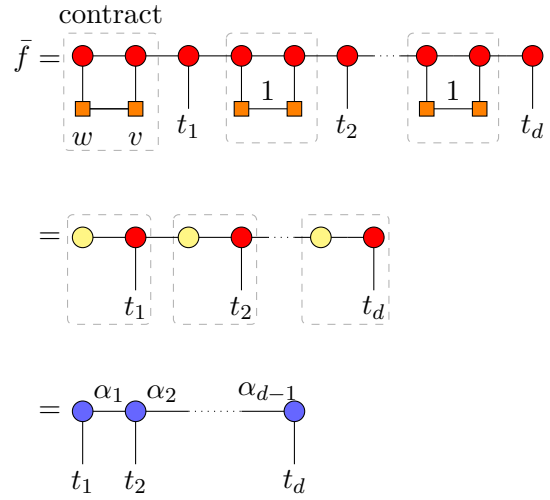
\begin{figure}[h!]
  \centering
  \begin{tikzpicture}

    \node[anchor=east] at (-0.15,0.0) {$\bar f=$};

    \draw[dashed,gray!60,rounded corners=2pt] (-0.25,0.30) rectangle (1.0,-1.25);
    \draw[dashed,gray!60,rounded corners=2pt] (1.9,0.30) rectangle (3.1,-1.15);
    \draw[dashed,gray!60,rounded corners=2pt] (4.3,0.30) rectangle (5.5,-1.15);

    \node[above] at (0.375,0.30) {contract};

    \draw (0,0) -- (0.7,0);
    \draw (0.7,0) -- (1.4,0);
    \draw (1.4,0) -- (2.1,0);
    \draw (2.1,0) -- (2.8,0);
    \draw (2.8,0) -- (3.5,0);
    \draw (3.5,0) -- (3.85,0);
    \draw[dotted] (3.85,0) -- (4.2,0);
    \draw (4.2,0) -- (4.9,0);
    \draw (4.9,0) -- (5.6,0);
    \draw (5.6,0) -- (5.95,0);

    \draw (0,-0.7) -- node[label={[shift={(-0.35,-0.75)}]$w$}]{} (0.7,-0.7);
    \draw (0,-0.7) -- node[label={[shift={(0.35,-0.75)}]$v$}]{} (0.7,-0.7);
    \draw (2.1,-0.7) -- node[label={[shift={(0,-0.15)}]$1$}]{} (2.8,-0.7);
    \draw (4.55,-0.7) -- node[label={[shift={(0,-0.15)}]$1$}]{} (5.25,-0.7);

    \draw (0,0) -- (0,-0.7);
    \draw (0.7,0) -- (0.7,-0.7);
    \draw (1.4,0) -- node[label={[shift={(0,-1)}]$t_1$}]{} ++(0,-0.7);
    \draw (2.1,0) -- (2.1,-0.7);
    \draw (2.8,0) -- (2.8,-0.7);
    \draw (3.5,0) -- node[label={[shift={(0,-1)}]$t_2$}]{} ++(0,-0.7);
    \draw (4.55,0) -- (4.55,-0.7);
    \draw (5.25,0) -- (5.25,-0.7);
    \draw (5.95,0) -- node[label={[shift={(0,-1)}]$t_d$}]{} ++(0,-0.7);

    \foreach \x in {0, 0.7, 1.4, 2.1, 2.8, 3.5, 4.55, 5.25, 5.95}
      \node[draw,circle,fill=red,scale=0.75] at (\x,0){};

    \foreach \x in {0, 0.7, 2.1, 2.8, 4.55, 5.25}
      \node[draw,rectangle,fill=orange,scale=0.75] at (\x,-0.7){};

    \begin{scope}[yshift=-2.4cm]

      \node[anchor=east] at (-0.15,0.0) {$=$};

      \draw[dashed,gray!60,rounded corners=2pt] (-0.2,0.30) rectangle (0.9,-1.15);
      \draw[dashed,gray!60,rounded corners=2pt] (1.2,0.30) rectangle (2.3,-1.15);
      \draw[dashed,gray!60,rounded corners=2pt] (2.9,0.30) rectangle (4.05,-1.15);

      \draw (0,0) -- (0.7,0);
      \draw (0.7,0) -- (1.4,0);
      \draw (1.4,0) -- (2.1,0);
      \draw (2.1,0) -- (2.45,0);
      \draw[dotted] (2.45,0) -- (2.8,0);
      \draw (2.8,0) -- (3.5,0);
      \draw (3.5,0) -- (3.85,0);

      \draw (0.7,0) -- node[label={[shift={(0,-1)}]$t_1$}]{} ++(0,-0.7);
      \draw (2.1,0) -- node[label={[shift={(0,-1)}]$t_2$}]{} ++(0,-0.7);
      \draw (3.85,0) -- node[label={[shift={(0,-1)}]$t_d$}]{} ++(0,-0.7);

      \foreach \x in {0, 1.4, 3.15}
        \node[draw,circle,fill=yellow!60,scale=0.75] at (\x,0){};
      \foreach \x in {0.7, 2.1, 3.85}
        \node[draw,circle,fill=red,scale=0.75] at (\x,0){};

      \node[anchor=east] at (-0.15,-2.3) {$=$};

      \draw (0,-2.3) -- node[label={[shift={(0,-0.15)}]$\alpha_1$}]{} (0.7,-2.3);
      \draw (0.7,-2.3) -- node[label={[shift={(0,-0.15)}]$\alpha_2$}]{} (1.4,-2.3);
      \draw[dotted] (1.4,-2.3) -- (2.1,-2.3);
      \draw (2.1,-2.3) -- node[label={[shift={(0,-0.15)}]$\alpha_{d-1}$}]{} (2.8,-2.3);

      \draw (0,-2.3) -- node[label={[shift={(0,-1)}]$t_1$}]{} ++(0,-0.7);
      \draw (0.7,-2.3) -- node[label={[shift={(0,-1)}]$t_2$}]{} ++(0,-0.7);
      \draw (2.8,-2.3) -- node[label={[shift={(0,-1)}]$t_d$}]{} ++(0,-0.7);

      \foreach \x in {0, 0.7, 2.8}
        \node[draw,circle,fill=blue!60,scale=0.75] at (\x,-2.3){};

    \end{scope}

  \end{tikzpicture}
  \caption{Contraction scheme for stage 2 of the algorithm. \textit{Top}: the $3d$ TT (red nodes) with the Euler weight vectors $w$ and $v$ (orange squares) attached to the $p_k$ and $j_k$ cores of each triplet. \textit{Middle}: after contracting $w$ and $v$ into each triplet, the two quadrature cores collapse into a single weight matrix (yellow nodes) that is multiplied into the $t_k$ core (red). \textit{Bottom}: the resulting $d$-dimensional time TT with bond dimensions $\alpha_1,\ldots,\alpha_{d-1}$. Because the weight vectors act independently on each triplet, no bond inflation occurs.}
  \label{fig:contraction_scheme}
\end{figure}

The weight vectors encode the Euler summation from~\eqref{eq:invmultidim}. The vector $w$ carries the complex phase factor and $v$ carries the Euler acceleration coefficients $c_j$ together with the alternating sign:
\begin{align}
    w[p] &= e^{-ip\pi/\ell}, \quad p = 1,\ldots,\ell,\\
    v[j] &= c_j\,(-1)^j, \quad j = -J,\ldots,J.
\end{align}
For each dimension $k$, we form the weighted sums over the physical indices of the $p$- and $j$-cores:
\begin{align}
    \hat{P}^{(k)} &= \sum_{p=1}^{\ell} w[p]\, P^{(k)}[:,p,:] \;\in\; \mathbb{C}^{r_{3k-3}\times r_{3k-2}},\\
    \hat{J}^{(k)} &= \sum_{j=-J}^{J} v[j]\, J^{(k)}[:,j,:] \;\in\; \mathbb{C}^{r_{3k-2}\times r_{3k-1}},
\end{align}
and multiply them to obtain the weight matrix
\begin{equation}\label{eq:wcore}
    W^{(k)} = \hat{P}^{(k)} \cdot \hat{J}^{(k)} \;\in\; \mathbb{C}^{r_{3k-3}\times r_{3k-1}}.
\end{equation}
We then absorb the time prefactor $c(t_i) = e^{A/(2\ell)}/(2t_i\ell)$ from~\eqref{eq:invmultidim}. For each grid point $t_i$, the matrix $T^{(k)}[:,t_i,:] \in \mathbb{C}^{r_{3k-1} \times r_{3k}}$ is the slice of the $t$-core at $t_i$. The contracted time core is
\begin{equation}\label{eq:time_core}
    B^{(k)}[:,t_i,:] = c(t_i)\; W^{(k)} \cdot T^{(k)}[:,t_i,:] \;\in\; \mathbb{C}^{r_{3k-3}\times r_{3k}}.
\end{equation}
The $d$ time cores $\{B^{(k)}\}_{k=1}^d$ form the final TT approximation of $\bar{f}(t_1,\ldots,t_d)$, with bond dimensions $r_3,r_6,\ldots,r_{3(d-1)}$. The contraction is summarized in Figure~\ref{fig:contraction_scheme}.

The weight vectors $w$ and $v$ are shared across all dimensions and couple only the quadrature cores within each triplet. Because the contraction is separable in the indices $p_k$ and $j_k$, it introduces no additional inter-dimensional correlations. Consequently, the inter-triplet bond dimensions of the TT representation of $\bar f$ are inherited from those of the TT representation of $\tilde f$. The complete contraction scheme is shown in Figure~\ref{fig:contraction_scheme}.

Contracting the quadrature legs and absorbing the prefactor requires $\mathcal{O}(d(\ell+J)\chi^2)$ operations to form the weight matrices $W^{(k)}$ in \eqref{eq:wcore}, and an additional $\mathcal{O}(d n_t \chi^2)$ operations to construct the time cores $B^{(k)}$. Once assembled, the resulting TT representation encodes $\bar{f}(t_1,\ldots,t_d)$ on the full evaluation grid, allowing statistical quantities to be extracted directly.

Under a low-rank assumption on $\tilde{f}$, the overall complexity remains inexpensive compared with the naive $\mathcal{O}((\ell J)^d)$ cost of direct summation. The dominant expense is the construction of the TT in Figure~\ref{fig:tt-structure}, which requires approximately $\mathcal{O}(d\,n_t\,\ell\,J\,\chi^2)$ evaluations of $\tilde{f}$, where $\chi$ denotes the maximum bond dimension. For example, with $d=5$, $\ell=10$, $J=25$, $n_t=20$, and $\chi=50$, TT-cross requires roughly $10^8$ evaluations of $\tilde{f}$, compared with $(\ell J)^d \approx 10^{12}$ evaluations for direct summation.

\begin{table}[t]
\centering
\small
\begin{tabular}{p{0.42\columnwidth}p{0.48\columnwidth}}
\hline
Operation & Leading cost \\
\hline
Direct tensor-product inversion & $\mathcal{O}((\ell J)^d)$ evaluations of $\tilde f$ \\
TT-cross construction of $\mathcal{F}$ & $\mathcal{O}(d\,n_t\,\ell\,N_j\,\chi^2)$ evaluations of $\tilde f$ \\
Quadrature contraction & $\mathcal{O}(d(\ell+J)\chi^2+d n_t\chi^2)$ operations \\
Storage of the final time TT & $\mathcal{O}(d n_t\chi^2)$ entries \\
\hline
\end{tabular}
\caption{Leading complexity contributions for the tensor-train inversion procedure, assuming a uniform time grid of size $n_t$, quadrature sizes $\ell$ and $J$, and maximum TT bond dimension $\chi$.}
\label{tab:complexity_summary}
\end{table}

\section{Error sources and parameter selection}\label{sec:error_bounds}
The total error of the inversion algorithm has four main contributions:
(i) the tensor-train approximation error introduced by TT-cross interpolation,
(ii) the discretization error from approximating the inversion contour,
(iii) the truncation error from replacing the infinite Euler series by a finite sum, and
(iv) floating-point round-off error. For the positive-support examples, we typically use $A=26.8$, giving discretization errors of order $\mathcal{O}(10^{-11})$. In the MNIG example, the bilateral admissibility condition \eqref{eq:MNIG_LT} restricts the allowable real contour shift, so we use $A=15$, yielding a discretization floor of order $\mathcal{O}(10^{-6})$.

\begin{enumerate}

\item \textbf{Cross interpolation error}

The cross-interpolation error arises from approximating $\tilde{f}(\vec{s})$ on the quadrature grid by a tensor train with finite bond dimension $\chi$. The TT-cross sweep is terminated when the relative error on a held-out validation set falls below $\varepsilon$. This gives an empirical estimate of the approximation error on the quadrature grid, but it does not imply a uniform error bound over all grid points. The worst-case error bound for the resulting TT approximation is $\mathcal{O}(\varepsilon\,\chi^d)$ \cite{qin2022error}, although in practice the algorithm performs substantially better when $\tilde{f}$ has low-rank structure on the grid, as in the Laplace transforms considered here. In all experiments we set $\varepsilon=10^{-6}$; higher accuracy is possible at the cost of larger bond dimension $\chi$.
    
\item \textbf{Discretization error}. 
The discretization error arises when the continuous inversion integral is approximated by a trapezoidal sum on an equally spaced frequency grid. As in \cite{MultidimensionalLaplace}, replacing the integral by a sum is justified through the Poisson summation formula. The trapezoidal approximation therefore recovers not only the desired value $f(t)$, but also a series of shifted aliasing terms. These terms are suppressed by the damping parameter $A_k$. In general, the discretization error satisfies
\begin{equation}\label{eq:discretization_error}
    \abs{\bar{e}} \leq C \sum_{k=1}^d \exp(-A_k),
\end{equation}
where $C \geq \abs{f}$ is a constant. When the same damping parameter $A$ is used in all dimensions, \eqref{eq:discretization_error} reduces to
\begin{equation}\label{eq:discretization}
    \abs{\bar{e}} \lessapprox Cd\exp({-A}).
\end{equation}
Thus increasing $A$ reduces the discretization error exponentially, while the dimensional contribution grows only linearly in $d$ for fixed $A$.

In our examples we used $A$ ranging from 15 to 26.8 which yields a discretization error of $\mathcal{O}(10^{-6})$ to $\mathcal{O}(10^{-11})$, respectively.

    \item \textbf{Infinite series approximation error}
For each dimension in \eqref{eq:invmultidim}, the infinite series is approximated by Euler summation, which introduces a truncation error controlled by the number of terms retained. As noted in \cite{abate1998numerical}, this effect is significant mainly near discontinuities in one or more dimensions. Following \cite{abate1992fourier}, we estimate the Euler-summation error by
\begin{equation}
    \varepsilon_{\text{trunc}}=\abs{E(m,n) - E(m-1, n)} \leq \frac{\Delta^m \alpha_{n+1}}{2^m},
\end{equation}
where $\Delta \alpha_j=\alpha_{j+1}-\alpha_j$ and $\Delta^k$ denotes $k$ applications of the forward-difference operator. Even for the modest choice $m=12$ and $n=5$, we obtain a truncation error of order $\mathcal{O}(10^{-11})$ at random sample points in the four-dimensional Wishart example shown in Figure~\ref{fig:Wishart4D}.

    \item \textbf{Round-off error} 
    
Round-off error is controlled by the Euler-summation parameter $\ell_k$. It arises because the prefactor $e^{A_k/(2\ell_k)}/(2\ell_k t_k)$ multiplies small, rapidly oscillating evaluations of $\tilde{f}$. This means significant digits can be lost when the prefactor is large. The damping parameter $A_k$ is fixed by the discretization-error target, so $\ell_k$ is the remaining control. Increasing $\ell_k$ shrinks the prefactor and reduces round-off error; the trade-off is that $\ell_k$ also determines the number of quadrature nodes per dimension, so larger $\ell_k$ increases $N$ and may increase the bond dimension required to represent $\tilde{f}$ on the finer quadrature grid; see \cite{scaling}.

\end{enumerate}

\section{Examples}\label{sec:examples}
We validate the tensor-train (TT) inversion algorithm introduced in Section~\ref{sec:tensor_formulation} by comparing recovered densities against Monte Carlo estimates and, where available, exact reference densities. The numerical experiments demonstrate that the method accurately reconstructs multivariate probability densities in dimensions that are beyond the reach of direct inversion. We consider three classes of examples: (i) multivariate normal-inverse Gaussian (MNIG) distributions, (ii) Wishart-type distributions with non-separable Laplace transforms, and (iii) correlated Gamma-type models.

To benchmark the TT inversion against a sampling-based approach, we draw independent samples from each target distribution and estimate the density using multivariate kernel density estimation (KDE). For the Wishart example, exact samples of the diagonal are obtained by sampling from the Wishart distribution directly and extracting the diagonal entries. For the MNIG distribution, we exploit its variance-mean mixture structure: a sample is generated by first drawing a mixing variable
\[
V \sim \mathrm{IG}(\delta/\gamma,\,\delta^2),
\]
and then setting
\begin{equation}
    X = \mu + (\Sigma\beta)\,V + \sqrt{V}\,L\,Z,
    \quad Z \sim \mathcal{N}(0,I),
\end{equation}
where $L$ is the Cholesky factor of $\Sigma$. This procedure yields exact independent samples at a computational cost of $\mathcal{O}(d)$ per draw.

Given $N$ samples $\{x^{(i)}\}_{i=1}^N \subset \mathbb{R}^d$, the density at a query point $t$ is estimated by placing a small Gaussian bump at each sample and averaging:
\begin{equation}
    \hat{f}(t) = \frac{1}{N}\sum_{i=1}^{N}
    \mathcal{N}\!\left(t;\, x^{(i)},\, H^2\right),
\end{equation}
where $\mathcal{N}(t;\mu,\Sigma)$ is the multivariate Gaussian density with mean $\mu$ and covariance $\Sigma$, and $H = \mathrm{diag}(h_1,\ldots,h_d)$ is a diagonal bandwidth matrix. The bandwidth in each coordinate is set by
Silverman's rule,
\begin{equation}
    h_k = N^{-1/(d+4)}\,\hat{\sigma}_k,
\end{equation}
where $\hat{\sigma}_k$ is the empirical standard deviation of the $k$-th coordinate across the $N$ samples. Larger $N$ shrinks the bandwidth and sharpens the estimate, but the rate of improvement slows dramatically in high dimensions: the root mean squared error (RMSE) of KDE scales as $\mathcal{O}(N^{-2/(d+4)})$ under standard smoothness assumptions. To quantify the variability of the Monte Carlo estimator, we ran each configuration 10 times using independent sample draws. The shaded band around each MC curve spans the pointwise minimum to maximum across all 10 repetitions, while the central line shows the mean.

\subsection{Multivariate normal-inverse Gaussian distribution}\label{sec:MNIG}
 
A multivariate normal-inverse Gaussian (MNIG) distribution is a variance-mean mixture of
a multivariate Gaussian distribution with an inverse Gaussian mixing distribution. Such
distributions have been applied in fields such as finance and signal processing
\cite{barndorff1994normal, hanssen2001normal, barndorff1982normal, oigaard2005estimation, oigard2002multivariate}.
 
An MNIG-distributed random vector $\vec{X}$ can be constructed as a variance-mean mixture
of a $d$-dimensional Gaussian random variable $\vec{Y}$ with a univariate inverse Gaussian
distributed mixing variable $Z$. We consider parameters $\alpha > 0$,
$\bm{\beta} \in \mathds{R}^d$, $\delta > 0$, $\bm{\mu} \in \mathds{R}^d$, and
$\Sigma \in \mathds{R}^{d \times d}$ symmetric positive definite, satisfying the
admissibility condition
\begin{equation}
    \alpha^2 > \bm{\beta}^\top \Sigma \bm{\beta}.
\end{equation}
The random vector $\vec{X}$ is defined as
\begin{equation}\label{eq:NIG}
    \vec{X} = \bm{\mu} + Z\Sigma\bm{\beta} + \sqrt{Z}\,\Sigma^{1/2}\vec{Y}
\end{equation}
where $Z \sim \mathrm{IG}(\delta^2,\, \alpha^2 - \bm{\beta}^\top \Sigma \bm{\beta})$ and
$\vec{Y} \sim \mathcal{N}_d(\vec{0}, I_d)$. The inverse Gaussian distribution
$\mathrm{IG}(\phi, \psi)$, $\phi, \psi > 0$, has density \cite{folks1978inverse}
\begin{equation}
    f_Z(z) = \left( \frac{\phi}{2\pi z^3} \right)^{1/2}
    \exp\!\left(\sqrt{\phi\psi} - \frac{1}{2}\!\left(\phi z^{-1} + \psi z\right)\right),
\end{equation}
for $z > 0$. From \eqref{eq:NIG}, the probability density function of $\vec{X}$ is given by
\begin{align}\label{eq:MNIG}
    f_{\vec{X}}(\vec{x})
    &= \int_0^\infty f_{\vec{X}|Z}(\vec{x}\mid z)\,f_Z(z)\,\mathrm{d}z \notag \\
    &= \frac{2\delta}{(2\pi)^{(d+1)/2}|\Sigma|^{1/2}}
    \left(\frac{\alpha}{Q(\vec{x})}\right)^{(d+1)/2} \notag \\
    & \hspace{4em} \times \exp\!\big(P(\vec{x})\big)\,
    K_{(d+1)/2}\,\big (\alpha  Q(\vec{x})\big),
\end{align}
where $K_\nu$ denotes the modified Bessel function of the second kind and
\begin{align}
    P(\vec{x}) &= \delta\sqrt{\alpha^2 - \bm{\beta}^\top \Sigma\, \bm{\beta}}
                 + \bm{\beta}^\top(\vec{x} - \bm{\mu}), \\
    Q(\vec{x}) &= \sqrt{\delta^2 + (\vec{x}-\bm{\mu})^\top \Sigma^{-1}(\vec{x}-\bm{\mu})}.
\end{align}
The bilateral Laplace transform of $\vec{X}$ is \cite{hanssen2001normal,
oigard2002multivariate}
\begin{align}\label{eq:MNIG_LT}
    \tilde{f}(\vec{s})
    & = \mathbb{E}\!\left[\exp\!\left(-\vec{s}^\top \vec{X}\right)\right] \notag \\
    & = \exp\!\Biggl(-\bm{\mu}^\top \vec{s}
        + \delta\Bigl(\sqrt{\alpha^2 - \bm{\beta}^\top \Sigma\,\bm{\beta}} \notag \\
    & \hspace{5em} - \sqrt{\alpha^2 - (\bm{\beta} - \vec{s})^\top \Sigma\,(\bm{\beta} - \vec{s})}\Bigr)\Biggr),
\end{align}
for $\vec{s} \in \mathbb{C}^d$ with $\vec{a} = \Re(\vec{s})$ satisfying
\begin{equation}\label{eq:MNIG_condition}
    (\bm{\beta} - \vec{a})^\top \Sigma\,(\bm{\beta} - \vec{a}) < \alpha^2.
\end{equation}
Since the MNIG distribution is supported on $\mathbb{R}^d$, $\tilde{f}$ is interpreted as a bilateral Laplace transform. This imposes constraints when evaluating on the complex quadrature nodes \eqref{eq:quad_nodes} with real part $\Re(\xi_k) = A/(2\ell t_k)$, independent of $j_k$ and $p_k$. For the MNIG Laplace transform \eqref{eq:MNIG_LT} to be well-defined at the node vector $\bm{\xi} = (\xi_1, \ldots, \xi_d)$, the admissibility condition \eqref{eq:MNIG_condition} must hold at each quadrature node. We therefore restrict the time grid to $t \geq t_{\mathrm{floor}}$, where $t_{\mathrm{floor}}$ is the smallest $t$ for which all quadrature nodes satisfy equation \eqref{eq:MNIG_condition}.

\begin{figure}[ht]
  \centering
  \includegraphics[width=0.5\textwidth]{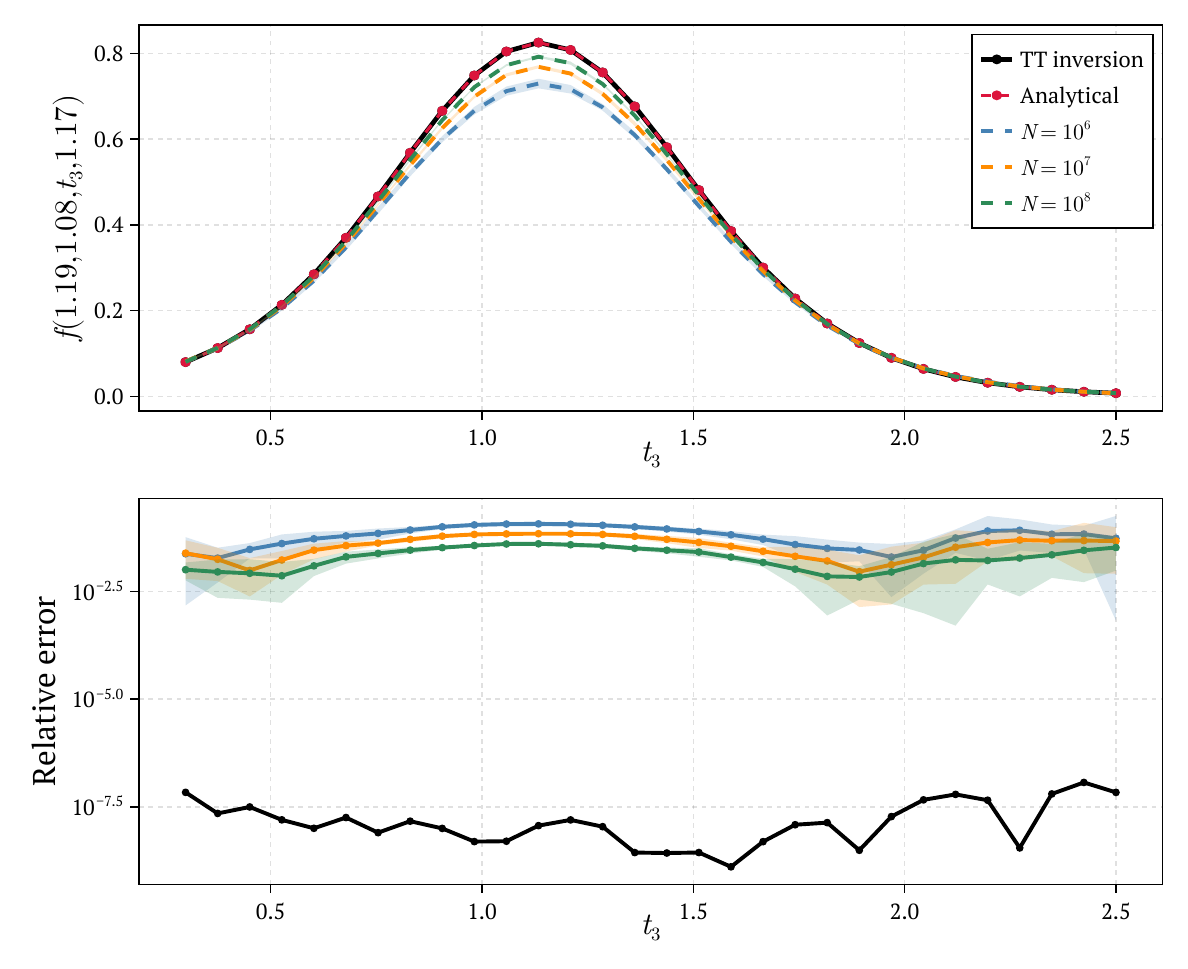}
  \caption{4-dimensional normal-inverse Gaussian distribution with a skewness parameter $\beta = [0.35, -0.2, 0.15, 0.25]$. The Monte Carlo samples are $N=10^6, 10^7, 10^8$ shown in blue, orange and green dashed lines, respectively.}
  \label{fig:MNIG_compare_MC}
\end{figure}

\begin{figure}[ht]
  \centering
  \includegraphics[width=0.5\textwidth]{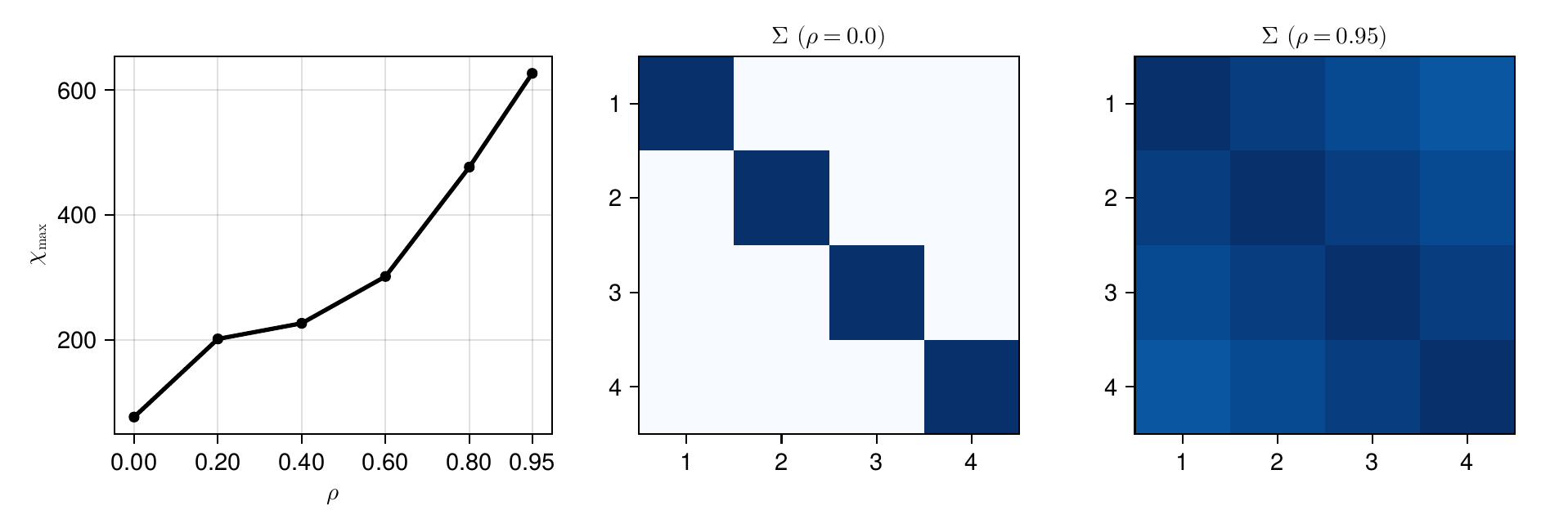}
  \caption{Maximum TT bond dimension $\chi$ as a function of the correlation parameter $\rho$ for a $d=4$ MNIG distribution (left), with the corresponding covariance matrix $\Sigma$ shown at $\rho=0$ (centre) and $\rho=0.95$ (right). To isolate the effect of correlation, $\alpha$ is adjusted at each $\rho$ so that $\gamma = \sqrt{\alpha^2 - \bm{\beta}^\top\Sigma\bm{\beta}}$ remains constant. The bond dimension grows substantially with $\rho$.}
  \label{fig:MNIG_correlationplot}
\end{figure}

\begin{figure}[ht]
  \centering
  \includegraphics[width=0.5\textwidth]{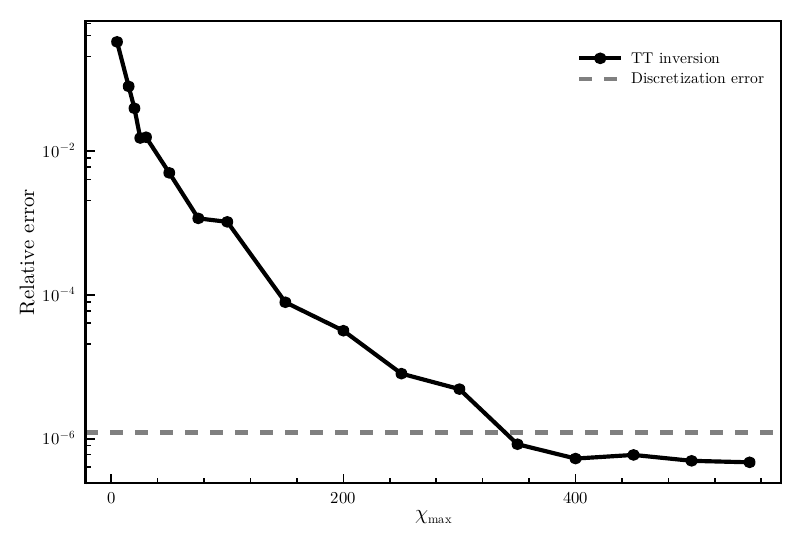}
  \caption{Relative error of the TT inversion versus maximum bond dimension $\chi$ for a $d=4$ MNIG distribution, with Laplace parameters $A=15$, $\ell=15$, $m=12$, $n_{\text{E}}=15$. The error decreases as the bond dimension grows, reaching the discretization floor (dashed) at $\chi \approx 450$, beyond which further increasing $\chi$ yields no improvement. The discretization floor $\approx 10^{-6}$ is set by the damping parameter $A$ and dimension $d$ via equation \eqref{eq:discretization}.}
  \label{fig:MNIG_convergenceplot}
\end{figure}

Figure~\ref{fig:MNIG_compare_MC} shows results for a four-dimensional normal-inverse Gaussian distribution with skewness parameter $\bm{\beta}=[0.35,-0.2,0.15,0.25]$. The density is evaluated on a one-dimensional slice along the third component $X_3$, with the remaining components fixed at their marginal means. Because the closed-form density \eqref{eq:MNIG} is available, we use it as the reference and report relative errors in the lower panel. We compare with Monte Carlo estimates using $N=10^6$, $10^7$, and $10^8$ samples, shown as blue, orange, and green dashed lines, respectively. The TT inversion closely matches the analytical density across the full range, whereas Monte Carlo requires $N=10^8$ samples to reach comparable accuracy.
To assess how correlation affects the required bond dimension, we use the covariance model $\Sigma_{ij}=\rho^{|i-j|}$. Figure~\ref{fig:MNIG_correlationplot} shows that as $\rho$ increases and the covariance matrix becomes less sparse, the maximum bond dimension $\chi$ must also increase. To isolate the effect of correlation, $\alpha$ is adjusted at each $\rho$ so that $\gamma=\sqrt{\alpha^2-\bm{\beta}^\top\Sigma\bm{\beta}}$ remains constant.
Figure~\ref{fig:MNIG_convergenceplot} shows the relative error as a function of the maximum bond dimension. We see that for a sufficiently high bond dimension we achieve a relative error lower than the discretization floor given by equation \eqref{eq:discretization}.

\begin{figure}[ht]
  \centering
  \includegraphics[width=0.5\textwidth]{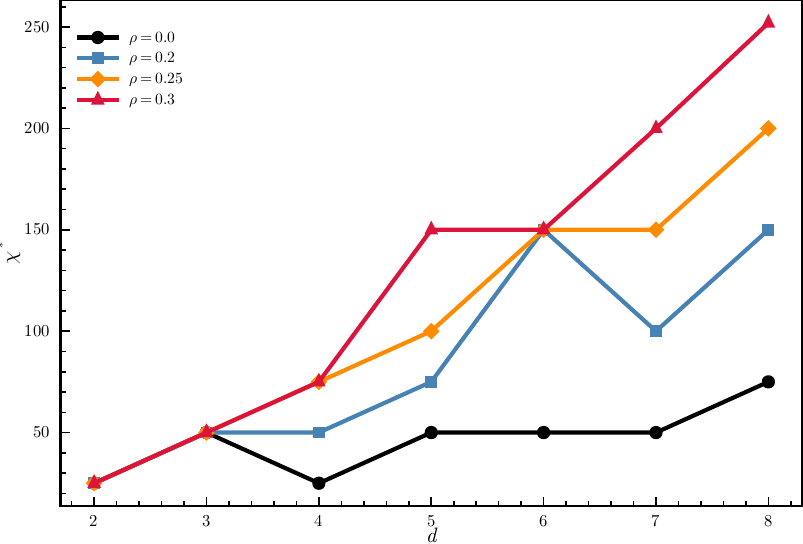}
  \caption{Minimum bond dimension $\chi$ required to achieve a relative error below $2 \times 10^{-3}$ as a function of dimension $d$, for the MNIG distribution with covariance at four correlation levels $\rho \in \{0.0, 0.2, 0.25, 0.3\}$. As in Figure~\ref{fig:MNIG_correlationplot}, $\alpha$ is adjusted at each $\rho$ to keep $\gamma = \sqrt{\alpha^2 - \bm{\beta}^\top\Sigma\bm{\beta}}$ fixed at $\gamma_0 = 3.9$. For $\rho = 0$ (no correlation) the required bond dimension grows slowly with $d$, reflecting that the uncorrelated MNIG is nearly separable. As correlation increases the required bond dimension grows faster, but remains moderate — reaching $\chi \approx 250$ at $d=8$, $\rho=0.3$.}
  \label{fig:MNIG_scalingplot}
\end{figure}

Figure~\ref{fig:MNIG_scalingplot} shows the minimum bond dimension required to obtain a relative error below $2\times10^{-3}$ as a function of the dimensionality, $d$ for different correlation parameters $\rho \in \{0.0, 0.2, 0.25, 0.3\}$. We see that as we increase the correlation we have to increase the bond dimension accordingly. The overlap at $d=6$ is because the bond dimension is incremented in discrete steps during \texttt{TT-cross}. A finer incrementation level would separate them at the cost of additional simulation time.

\begin{figure}[ht]
  \centering
  \includegraphics[width=0.5\textwidth]{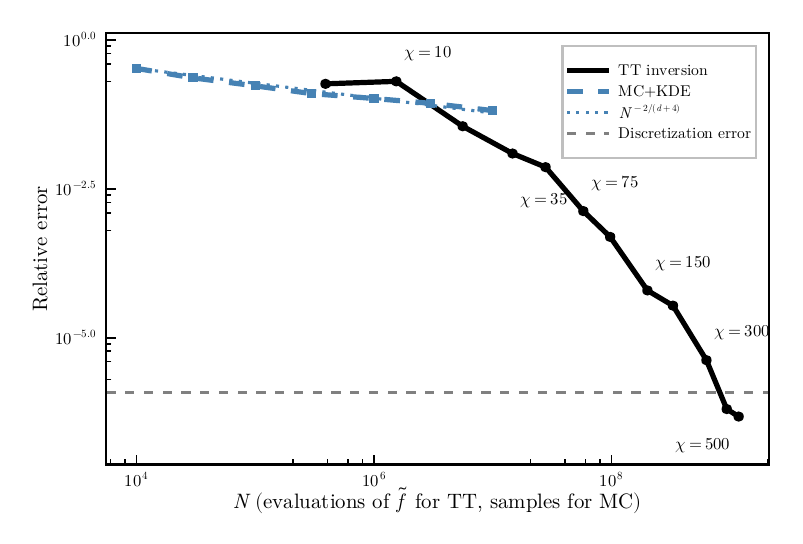}
  \caption{Cost-accuracy comparison between TT inversion and Monte Carlo KDE for the $d=4$ MNIG distribution. The horizontal axis shows the number of evaluations of $\tilde{f}$ (for TT) or samples drawn (for MC+KDE). The TT inversion (black) rapidly reaches the discretization floor $\approx 10^{-6}$ as the bond dimension $\chi$ grows, while MC+KDE (blue dashed) converges at the rate $N^{-2/(d+4)}$ (dotted), which is $N^{-1/4}$ for $d=4$. The TT inversion achieves several orders of magnitude lower error at the same evaluation budget, with the annotated $\chi$ values indicating the bond dimension at each point along the TT curve.}
  \label{fig:MNIG4D_cost}
\end{figure}

Figure~\ref{fig:MNIG4D_cost} shows a cost-accuracy comparison between the TT inversion algorithm and the Monte Carlo KDE for $d=4$. The horizontal axis shows the number of function evaluations of $\tilde{f}$ used by the \texttt{TT-cross} algorithm and the number of drawn samples for the MC+KDE method. 

The TT inversion algorithm obtains a relative error lower than the discretization error \eqref{eq:discretization} at a maximum bond dimension $\chi \approx 450$.

\subsection{The Wishart distribution}\label{sec:Wishart}

As a benchmark for non-separable multivariate Laplace transforms, we consider
the real Wishart distribution; see \cite{LETAC20081393}. Let
\begin{equation}
    S = XX^\top,
\end{equation}
where $X \in \mathbb{R}^{d \times \nu}$ with $\nu \ge d$, whose columns are
independent and identically distributed as $\mathcal{N}(0,\Sigma)$, and
$\Sigma \in \mathbb{R}^{d \times d}$ is a symmetric positive definite matrix.
Then $S \sim \mathcal{W}_d(\nu,\Sigma)$ and the Laplace transform of $S$ is given by
\begin{align}
    \phi(\theta)
    &= \mathbb{E}\!\left[\exp\!\left(\operatorname{tr}(\theta S)\right)\right] \\
    &= \mathbb{E}\!\left[\exp\!\left(\operatorname{tr}(X^\top \theta X)\right)\right],
\end{align}
for matrices $\theta \in \mathbb{R}^{d \times d}$ such that $I - 2\Sigma\theta$
is positive definite. Writing $X = \Sigma^{1/2} Y$ with $Y$ having i.i.d.\
standard normal entries, we obtain
\begin{equation}
    \phi(\theta) = \mathbb{E}\!\left[\exp\!\left(\operatorname{tr}(Y^\top A Y)
    \right)\right],
    \quad A = \Sigma^{1/2}\theta \Sigma^{1/2}.
\end{equation}
Using the spectral decomposition of $A$ and independence of the columns of $Y$,
this yields
\begin{equation}
    \phi(\theta) = \det(I - 2\Sigma \theta)^{-\nu/2}.
\end{equation}

Restricting to diagonal test matrices $\theta = -\operatorname{diag}(s)$, we
obtain the Laplace transform of the diagonal vector $(S_{11}, \dots, S_{dd})$:
\begin{equation}
\label{eq:Wishart_Laplace}
    \phi(s) = \det(I + 2\Sigma \operatorname{diag}(s))^{-\nu/2},
    \quad \Re(s_i) \ge 0.
\end{equation}
The transform \eqref{eq:Wishart_Laplace} is non-separable in the variables $s_i$ whenever $\Sigma$ has non-zero off-diagonal entries, because the determinant couples all variables simultaneously. Although the full Wishart density has a closed-form expression, the marginal density of $(S_{11},\dots,S_{dd})$ does not have a simple closed form when $\Sigma$ is non-diagonal. Our main comparison in this case is therefore based on Monte Carlo sampling.

To quantify how the coupling in $\Sigma$ affects the complexity of the TT
representation, we parametrize $\Sigma_{ij} = \rho^{|i-j|}$. For $\rho = 0$, $\Sigma$ is diagonal and \eqref{eq:Wishart_Laplace} factors into
a product of scalar transforms yielding a TT with small bond dimensions. As $\rho$ increases, off-diagonal entries grow and all variables become coupled through the determinant, so higher bond dimension is
required. This parametrization is used below as a controlled way to interpolate between the separable diagonal case and increasingly non-separable transforms.

\begin{figure}[ht]
  \centering
  \includegraphics[width=0.5\textwidth]{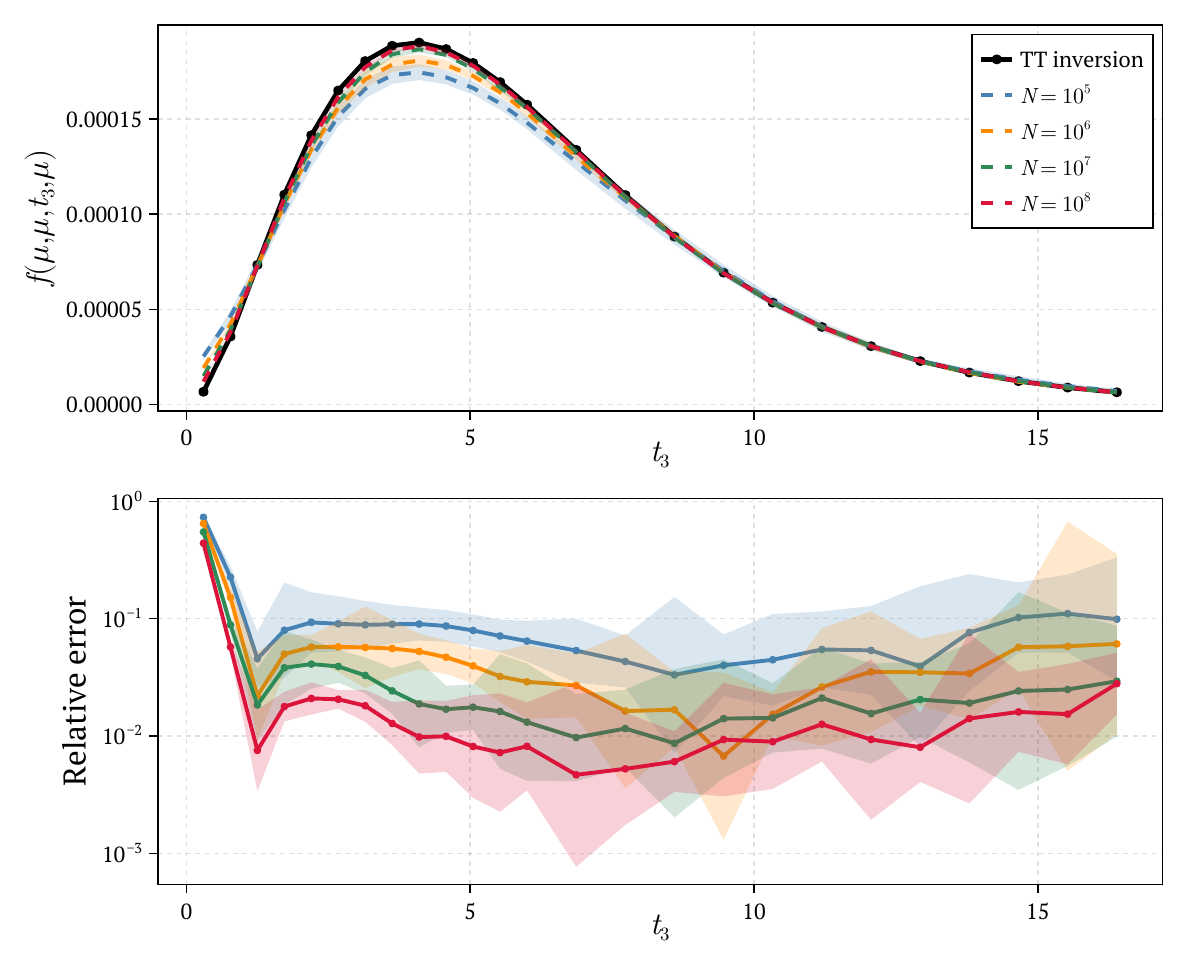}
  \caption{TT inversion applied to the $d = 4$ Wishart Laplace transform
    \eqref{eq:Wishart_Laplace} with $\nu = d+2$ and correlation parameter
    $\rho = 0.15$. The density is evaluated along the slice $t_k = \mu$ for
    $k \ne 3$, where $\mu = \nu\Sigma_{11} = 6$. Monte Carlo KDE estimates for
    $N = 10^5, 10^6, 10^7$ and $10^8$ are shown in dashed blue, orange, green
    and red, respectively. Parameters: $A = 26.8$, $\ell = 12$, $m = n = 10$. The lower panel shows the relative difference between KDE estimates and the TT inversion, not the true numerical error.}
  \label{fig:Wishart4D}
\end{figure}

We validate the inversion algorithm for a fixed correlation $\rho = 0.15$ with
$d = 4$ and $\nu = d+2$. Since no closed-form expression is available for the
joint density of the diagonal entries in this case, we compare against kernel
density estimates (KDE) obtained by drawing samples $S \sim \mathcal{W}_d(\nu, \Sigma)$
and extracting the diagonals. The density is evaluated along the one-dimensional
slice
\[
(t_1, \dots, t_d) = (\mu, \dots, \mu_{k-1}, t, \mu_{k+1}, \dots, \mu),
\]
where $\mu_i = \mathbb{E}[S_{ii}] = \nu\Sigma_{ii}$. Figure~\ref{fig:Wishart4D}
shows close agreement between the TT inversion and the KDE across all sample
sizes, with the KDE converging towards the TT result as $N$ increases.

\begin{figure}[ht]
  \centering
  \includegraphics[width=0.5\textwidth]{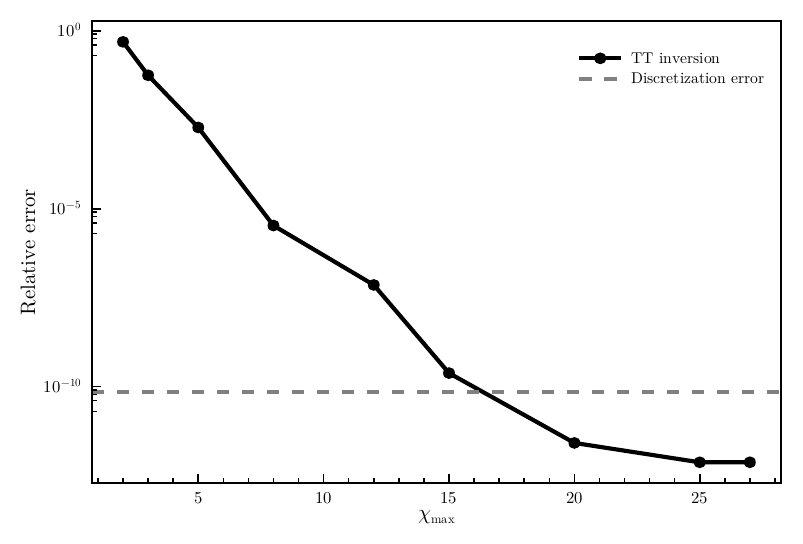}
  \caption{Relative error of the TT inversion versus maximum bond dimension $\chi$ for a $d=5$ Wishart distribution with $\Sigma = \operatorname{diag}(\lambda_1,\dots,\lambda_d)$ and inversion parameters $A=25$, $\ell=12$, $m=12$, $n_{\text{E}}=12$. The discretization error $\approx 10^{-10}$ is set by the damping parameter $A$ and dimension $d$ via equation \eqref{eq:discretization}.}
  \label{fig:Wishart_convergence}
\end{figure}

As a separate numerical check, we also considered the separable case
$\Sigma = \operatorname{diag}(\lambda_1,\dots,\lambda_d)$, for which the diagonal entries of $S$ are independent and the exact joint density is a product of Gamma distributions, $S_{kk}\sim \mathrm{Gamma}(\nu/2,2\lambda_k)$. In this case the relative error decreases rapidly with the bond dimension $\chi$ and saturates at
the discretization floor determined by the quadrature parameters, confirming that the dominant error source at convergence is the discretization error rather than the TT approximation itself. This is shown in Figure~\ref{fig:Wishart_convergence} and described in Section~\ref{sec:error_bounds}.

\subsection{Multivariate Correlated Gamma-Type Model}
We now consider a model motivated by factor models for portfolio credit risk, in particular the CreditRisk+ family of models introduced in \cite{suisse1997}. In such models, default counts or losses are often conditionally independent given a collection of latent sector or macroeconomic risk factors. Gamma-distributed factors are a natural choice because Gamma mixing of Poisson intensities leads to analytically tractable overdispersion and positive dependence across obligors or portfolio segments. A small number of common Gamma factors can therefore generate correlated losses across many names while retaining a simple transform representation.

The correlated Gamma-type model represents a $d$-dimensional non-negative random vector
$\boldsymbol{\Lambda} = (\Lambda_1,\dots,\Lambda_d)^\top$ as a linear mixture of $K$
independent Gamma factors. Let $Y_1,\dots,Y_K$ be independent with
$Y_k \sim \mathrm{Gamma}(\alpha_k,\beta_k)$, where $\alpha_k > 0$ is the shape
and $\beta_k > 0$ the rate, so that
\begin{equation}
  f_{Y_k}(y)
  = \frac{\beta_k^{\alpha_k}}{\Gamma(\alpha_k)}\,y^{\alpha_k-1}e^{-\beta_k y},
  \qquad y > 0.
\end{equation}
Given a $d\times K$ loading matrix $W = (w_{jk})\in\mathbb{R}_{\ge 0}^{d\times K}$,
each component is defined by
\begin{equation}
  \Lambda_j := \sum_{k=1}^K w_{jk} Y_k, \qquad j = 1,\dots,d,
\end{equation}
or in vector form $\boldsymbol{\Lambda} = W\mathbf{Y}$. Components that share
factor loadings are positively correlated, with the entire dependence structure
controlled by $W$ and the factor variances $\alpha_k/\beta_k^2$.

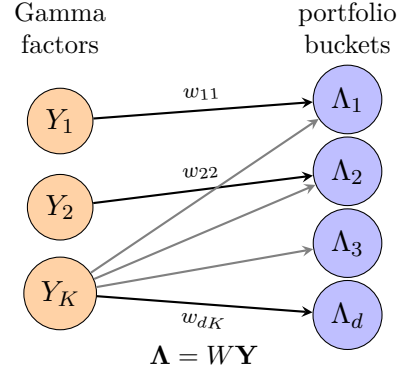
\begin{figure}[ht]
\centering
\begin{tikzpicture}[>=stealth,scale=0.95]
  \node[draw,circle,fill=orange!35,minimum size=0.55cm] (y1) at (0,1.2) {$Y_1$};
  \node[draw,circle,fill=orange!35,minimum size=0.55cm] (y2) at (0,0.0) {$Y_2$};
  \node[draw,circle,fill=orange!35,minimum size=0.55cm] (yc) at (0,-1.2) {$Y_K$};
  \node[draw,circle,fill=blue!25,minimum size=0.65cm] (l1) at (4,1.5) {$\Lambda_1$};
  \node[draw,circle,fill=blue!25,minimum size=0.65cm] (l2) at (4,0.5) {$\Lambda_2$};
  \node[draw,circle,fill=blue!25,minimum size=0.65cm] (l3) at (4,-0.5) {$\Lambda_3$};
  \node[draw,circle,fill=blue!25,minimum size=0.65cm] (ld) at (4,-1.5) {$\Lambda_d$};
  \draw[->,thick] (y1) -- node[above,sloped,font=\scriptsize] {$w_{11}$} (l1);
  \draw[->,thick] (y2) -- node[above,sloped,font=\scriptsize] {$w_{22}$} (l2);
  \draw[->,thick] (yc) -- node[below,sloped,font=\scriptsize] {$w_{dK}$} (ld);
  \draw[->,gray,thick] (yc) -- (l1);
  \draw[->,gray,thick] (yc) -- (l2);
  \draw[->,gray,thick] (yc) -- (l3);
  \node[align=center,font=\small] at (0,2.50) {Gamma\\factors};
  \node[align=center,font=\small] at (4,2.50) {portfolio\\buckets};
  \node[font=\small] at (2,-2.05) {$\boldsymbol{\Lambda}=W\mathbf{Y}$};
\end{tikzpicture}
\caption{Schematic of the correlated Gamma factor model. Idiosyncratic and common Gamma factors $Y_k$ feed portfolio buckets through the loading matrix $W$, producing correlated stochastic intensities $\Lambda_j$.}
\label{fig:gamma_factor_loading}
\end{figure}

In this interpretation, the factors $Y_k$ represent independent sector-level or systematic risk shocks, while the loading matrix $W$ maps these shocks into portfolio buckets, rating classes, regions, or other aggregated exposure groups as shown in Figure~\ref{fig:gamma_factor_loading}. The component $\Lambda_j$ may be interpreted as the stochastic intensity or loss driver for bucket $j$. Components with overlapping factor loadings are positively correlated, reflecting shared exposure to the same systematic credit-risk factors. This makes the model a useful test case for multidimensional Laplace inversion: the transform is simple and factorized, but the induced joint density of the correlated intensity vector is generally not available in closed form, especially when the number of latent factors exceeds the number of observed components.

This construction captures a central feature of CreditRisk+-type models: dependence is not introduced by specifying a dense covariance matrix directly, but by exposing several portfolio components to the same latent Gamma factors. The loadings $w_{jk}$ therefore have a direct risk interpretation, measuring the sensitivity of bucket $j$ to factor $k$. Tail events in several buckets are then driven by large realizations of common factors, making joint and conditional tail probabilities especially relevant diagnostics.

A useful feature of this construction is that the multivariate Laplace transform
$\phi_{\boldsymbol{\Lambda}}(\mathbf{s})
= \mathbb{E}[e^{-\mathbf{s}^\top\boldsymbol{\Lambda}}]$
factors over the $K$ independent components. Substituting
$\Lambda_j = \sum_k w_{jk}Y_k$ and swapping the order of summation gives
$\mathbf{s}^\top\boldsymbol{\Lambda} = \sum_{k=1}^K (\sum_j w_{jk}s_j)\,Y_k$,
and independence then yields
\begin{equation}
  \phi_{\boldsymbol{\Lambda}}(\mathbf{s})
  = \prod_{k=1}^K \mathbb{E}\!\left[
      e^{-\bigl(\sum_{j=1}^d w_{jk}s_j\bigr)Y_k}
    \right].
\end{equation}
Since $\mathbb{E}[e^{-tY_k}] = (\beta_k/(\beta_k+t))^{\alpha_k}$
for all $t > -\beta_k$, we obtain
\begin{equation}
  \label{eq:GammaTypeLaplace}
  \phi_{\boldsymbol{\Lambda}}(\mathbf{s})
  = \prod_{k=1}^K
    \left(
      \frac{\beta_k}{\beta_k + \sum_{j=1}^d w_{jk} s_j}
    \right)^{\!\alpha_k}\!,
\end{equation}
which is well-defined for all $\mathbf{s}\ge 0$ since $w_{jk}\ge 0$ and $\beta_k>0$,
and extends analytically to $\mathbf{s}\in\mathbb{C}^d$ with
$\beta_k + \sum_j w_{jk}\mathrm{Re}(s_j)>0$ for each $k$.
Evaluating~\eqref{eq:GammaTypeLaplace} at any quadrature node
costs $\mathcal{O}(Kd)$ operations, independent of whether $K\le d$ or $K>d$.

\begin{figure}[ht]
  \centering
  \includegraphics[width=0.5\textwidth]{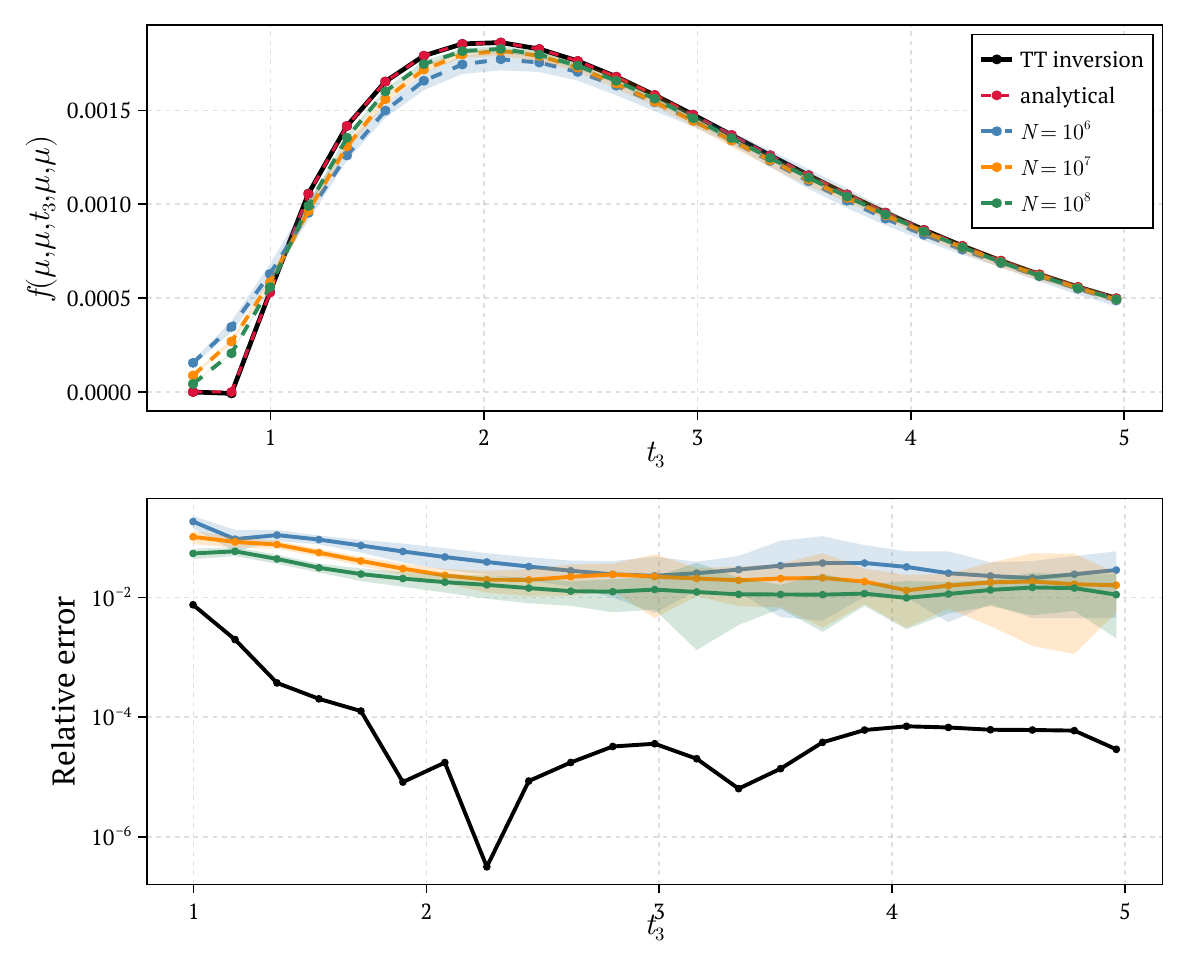}
  \caption{Correlated Gamma-type model, $d = K = 5$, with factors
    $Y_k \sim \mathrm{Gamma}(2,1)$ and loading matrix
    $W = (1-\varepsilon)I_5 + \varepsilon\,\mathbf{1}\mathbf{1}^\top$,
    $\varepsilon = 0.1$.
    The TT inversion (solid black, with inversion parameters $A=20, \ell = 12,, n_E=12,m_E=12$ with $\chi = 500$) is compared against the analytical density
    (dashed crimson, obtained via change of variables) and KDE with
    $N = 10^6$ (blue), $10^7$ (orange), and $10^8$ (green) samples,
    evaluated along the $\Lambda_3$ axis at the marginal means of the
    remaining components.}
  \label{fig:GammaType}
\end{figure}

We consider two instances. In the first, $d = K = 5$ with all factors
$Y_k\sim\mathrm{Gamma}(2,1)$ and loading matrix
$W = (1-\varepsilon)I_5 + \varepsilon\,\mathbf{1}\mathbf{1}^\top$, $\varepsilon = 0.1$:
each component draws primarily from its own factor with a 10\% cross-contribution
from the remaining four. Since $W$ is square and invertible, the joint density
admits the closed-form expression
\begin{equation}
  f_{\boldsymbol{\Lambda}}(\boldsymbol{\lambda})
  = \frac{f_{\mathbf{Y}}(W^{-1}\boldsymbol{\lambda})}{|\det W|},
\end{equation}
which serves as an exact validation reference.
Figure~\ref{fig:GammaType} shows the density along the $\Lambda_3$ axis with all
other components held at their marginal means. The TT inversion closely matches the
analytical reference, while the KDE estimator with $N=10^8$ samples still shows a
visible discrepancy near the mode — a consequence of the
$\mathcal{O}(N^{-2/9})$ pointwise convergence rate for $d=5$.

\begin{figure}[ht]
  \centering
  \includegraphics[width=0.5\textwidth]{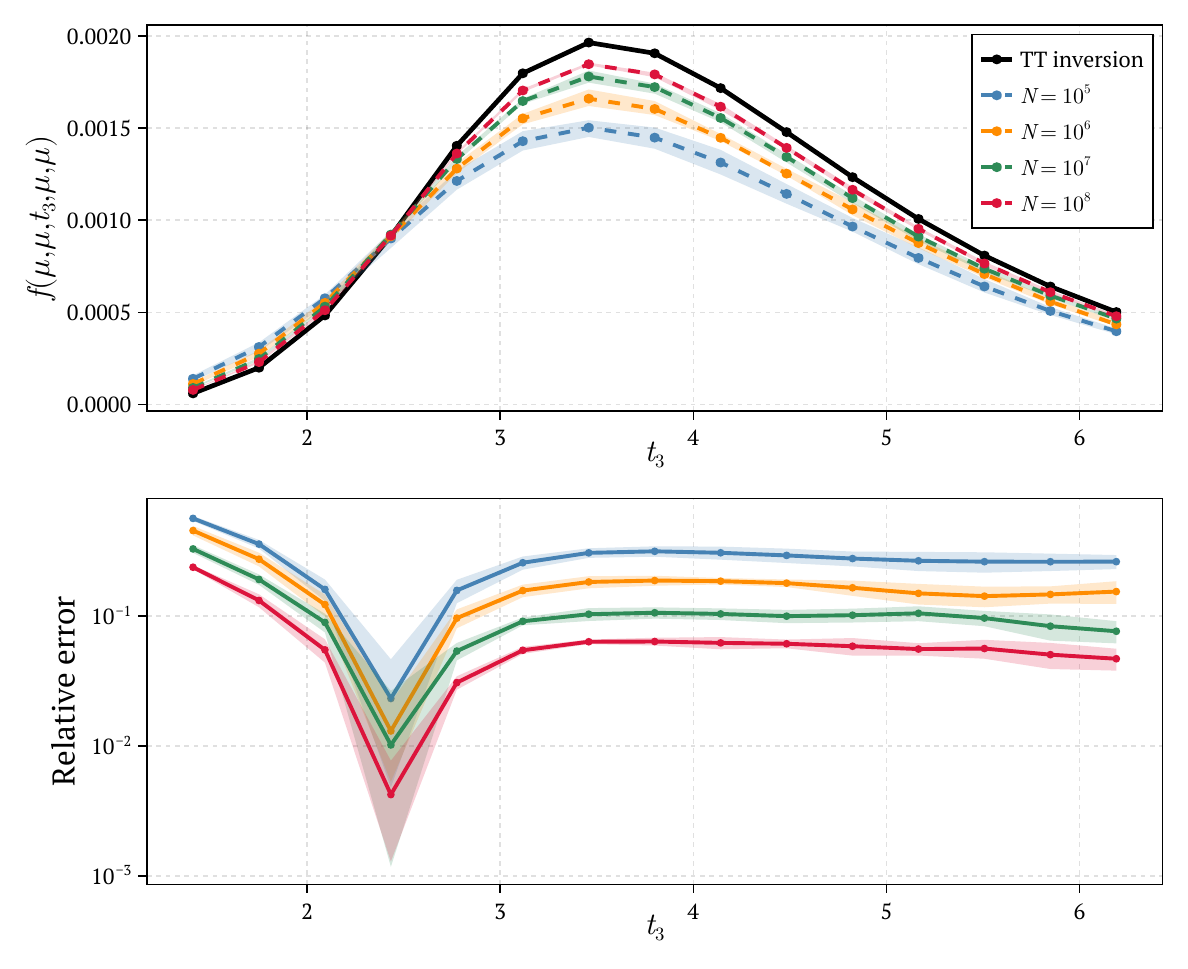}
  \caption{Correlated Gamma-type model, $d = 5$, $K = 8$, with factors
    $Y_k \sim \mathrm{Gamma}(2,1)$ and loading matrix
    $W = [I_{5\times 5}\;\mid\; 0.3\cdot\mathbf{1}_{5\times 3}]$
    (five idiosyncratic plus three common factors).
    No closed-form density exists for this  case.
    The TT inversion (solid black, inversion parameters $A=25,\ell=5,n_E=10,m_e=12)$ with $\chi = 250$) is compared against KDE with
    $N = 10^5$ (blue), $10^6$ (orange), $10^7$ (green), and $10^8$ (red) samples,
    evaluated along the $\Lambda_3$ axis at the marginal means of the
    remaining components. The lower panel shows the relative difference between KDE estimates and the TT inversion, not the true numerical error.}
  \label{fig:CorrGammad5k8}
\end{figure}

In the second instance we take $d = 5$ and $K = 8$, with loading matrix
$W = [I_{5\times 5}\;\mid\;0.3\cdot\mathbf{1}_{5\times 3}]$: five idiosyncratic
factors plus three common factors with loading $0.3$ shared across all components.
For $K > d$ no change-of-variables density formula is available, so the joint density
has no closed form. The TT inversion of~\eqref{eq:GammaTypeLaplace} therefore provides a way to evaluate $f_{\boldsymbol{\Lambda}}$ without relying on sampling.
Figure~\ref{fig:CorrGammad5k8} shows the result alongside KDE curves for
$N = 10^5$, $10^6$, $10^7$, and $10^8$ samples; the shared common factors induce
substantial off-diagonal covariance that slows KDE convergence, whereas the TT
inversion provides a stable density estimate at a fraction of the sampling cost.

\section{Distributional queries from the tensor-train density}\label{sec:distributional_queries}

The output of the inversion procedure is not only a pointwise approximation of the density. It is a reusable tensor-train representation of the joint distribution on the finite evaluation grid. Once this representation has been constructed, marginal densities, conditional distribution functions, and dependence measures can be computed by deterministic tensor contractions, without additional sampling or explicit construction of the full tensor-product array.

Since the inverse Laplace transform is evaluated on a finite tensor-product grid, the recovered tensor-train density is first normalized on this grid. We compute
\begin{equation}\label{eq:normalization}
Z =
\prod_{k=1}^{d} \Delta t_k
\sum_{i_1,\ldots,i_d}
f_{\mathrm{TT}}(i_1,\ldots,i_d),
\end{equation}
where $\Delta t_k$ is the grid spacing in dimension $k$. We use $f_{\mathrm{TT}}/Z$ in all subsequent distributional queries.

\subsection{Marginal densities}

Low-dimensional marginals are obtained by contracting over the unwanted tensor-train modes with the corresponding quadrature weights. For example, the two-dimensional marginal density of $(X_i,X_j)$ is approximated by
\begin{equation}\label{eq:marginals}
f_{ij}(t_i,t_j)
\approx
\frac{1}{Z}
\prod_{k\neq i,j}\Delta t_k
\sum_{\{i_k:k\neq i,j\}}
f_{\mathrm{TT}}(i_1,\ldots,i_d),
\end{equation}
where the indices corresponding to dimensions $i$ and $j$ are left uncontracted. One-dimensional marginal densities are obtained analogously by leaving a single physical index open.

\begin{figure}[ht]
\centering
\includegraphics[width=0.5\textwidth]{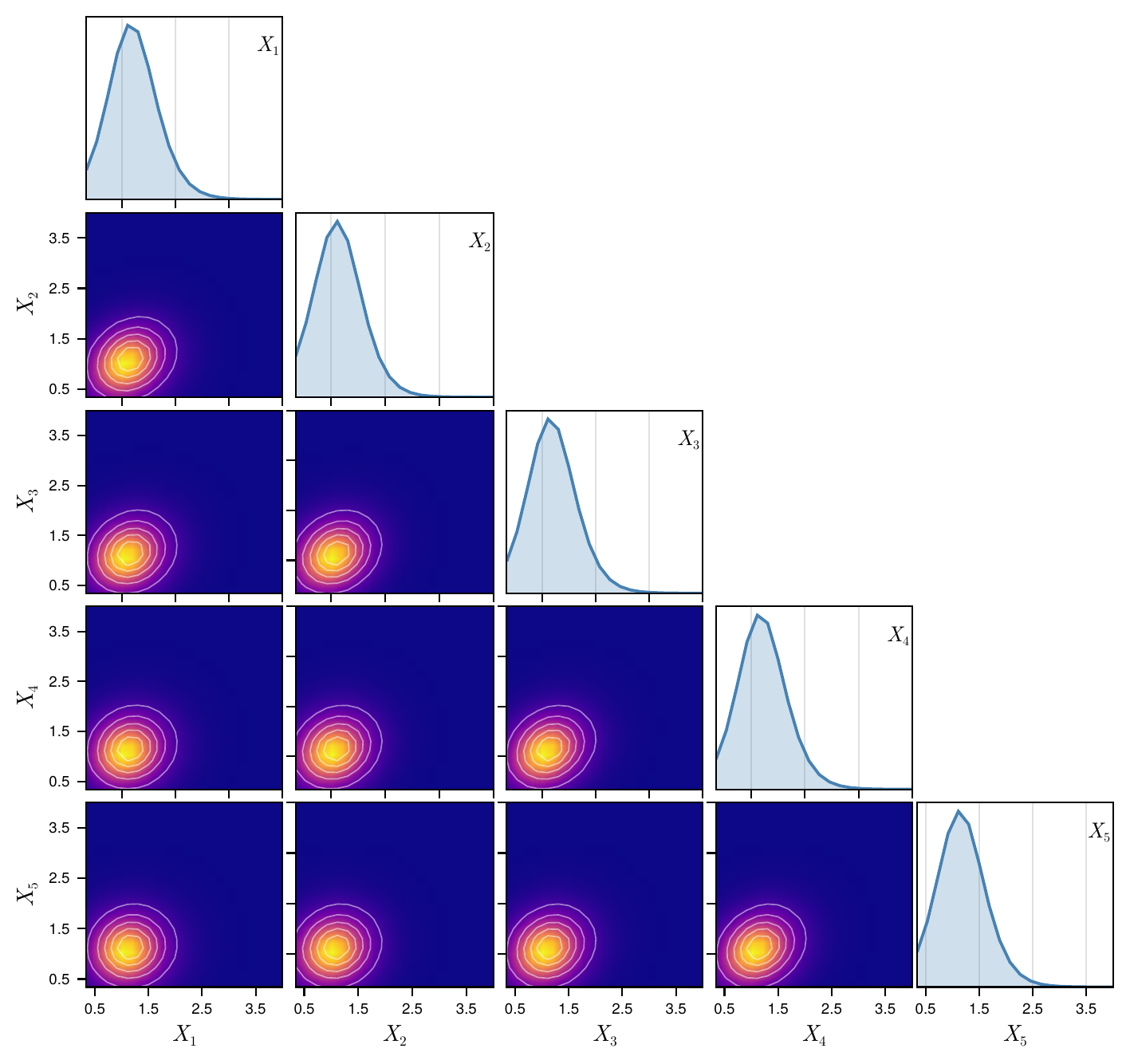}
\caption{Corner plot of the $d=5$ MNIG distribution recovered by TT inversion. Diagonal panels show the one-dimensional marginal density of each component $X_k$. Off-diagonal panels in the lower triangle show the two-dimensional joint marginals of $(X_j,X_i)$ as heat maps with contour lines, illustrating the pairwise dependence structure across all ten pairs of components. All marginals are obtained from the same normalized TT by dimension-wise contraction, with no additional sampling required. The parameters used are $\alpha=4.0$, $\delta=1.0$, and $\beta=[0.35,-0.2,0.15,0.25,0.10]$.}
\label{fig:marginals}
\end{figure}

Figure~\ref{fig:marginals} shows the resulting one- and two-dimensional marginals for the five-dimensional MNIG example. The diagonal panels show the marginal density of each component, while the off-diagonal panels show pairwise joint marginals. All panels are extracted from the same normalized TT representation, with no additional sampling and no explicit construction of the full $d$-dimensional array.

\subsection{Conditional cumulative distribution functions}

The TT representation also enables efficient computation of conditional
probabilities via marginalisation. Consider the joint conditional event
$\{X_2 > c_2,\, X_3 > c_3\}$; the conditional CDF of $X_1$ given this event is
\begin{multline}\label{eq:conditional}
F_{X_1\mid X_2>c_2,\,X_3>c_3}(x_1) = \mathbb{P}\!\left( X_1\leq x_1 \mid X_2>c_2,\, X_3>c_3\right) \\
= \frac{\mathbb{P}\!\left(X_1\leq x_1,\,
X_2>c_2,\,X_3>c_3\right)}{\mathbb{P}\!\left(X_2>c_2,\,X_3>c_3\right)}.
\end{multline}
Both numerator and denominator are obtained by contracting the TT cores with quadrature-weighted indicator vectors. Specifically, the numerator is computed as
\begin{equation}
  \int_{-\infty}^{x_1}\int_{c_2}^{\infty}\int_{c_3}^{\infty}
  f_{X_1,X_2,X_3}(u_1,u_2,u_3)\,\mathrm{d}u_3\,\mathrm{d}u_2\,\mathrm{d}u_1,
\end{equation}
where $f_{X_1,X_2,X_3}$ is the three-dimensional marginal density, obtained by
integrating the full five-dimensional TT density over $X_4$ and $X_5$.

We evaluate \eqref{eq:conditional} for the five-dimensional MNIG model of
Section~\ref{sec:MNIG}, with inversion parameters $A=15$, $\ell=15$, $m=15$,
$n_{\mathrm{E}}=15$, and thresholds $c_2=2.0$, $c_3=3.0$. The joint
conditional event $\{X_2>2.0,\,X_3>3.0\}$ is rare under the marginal
distribution, making it a challenging target for Monte Carlo estimation.

\begin{figure}[ht]
\centering
\includegraphics[width=0.5\textwidth]{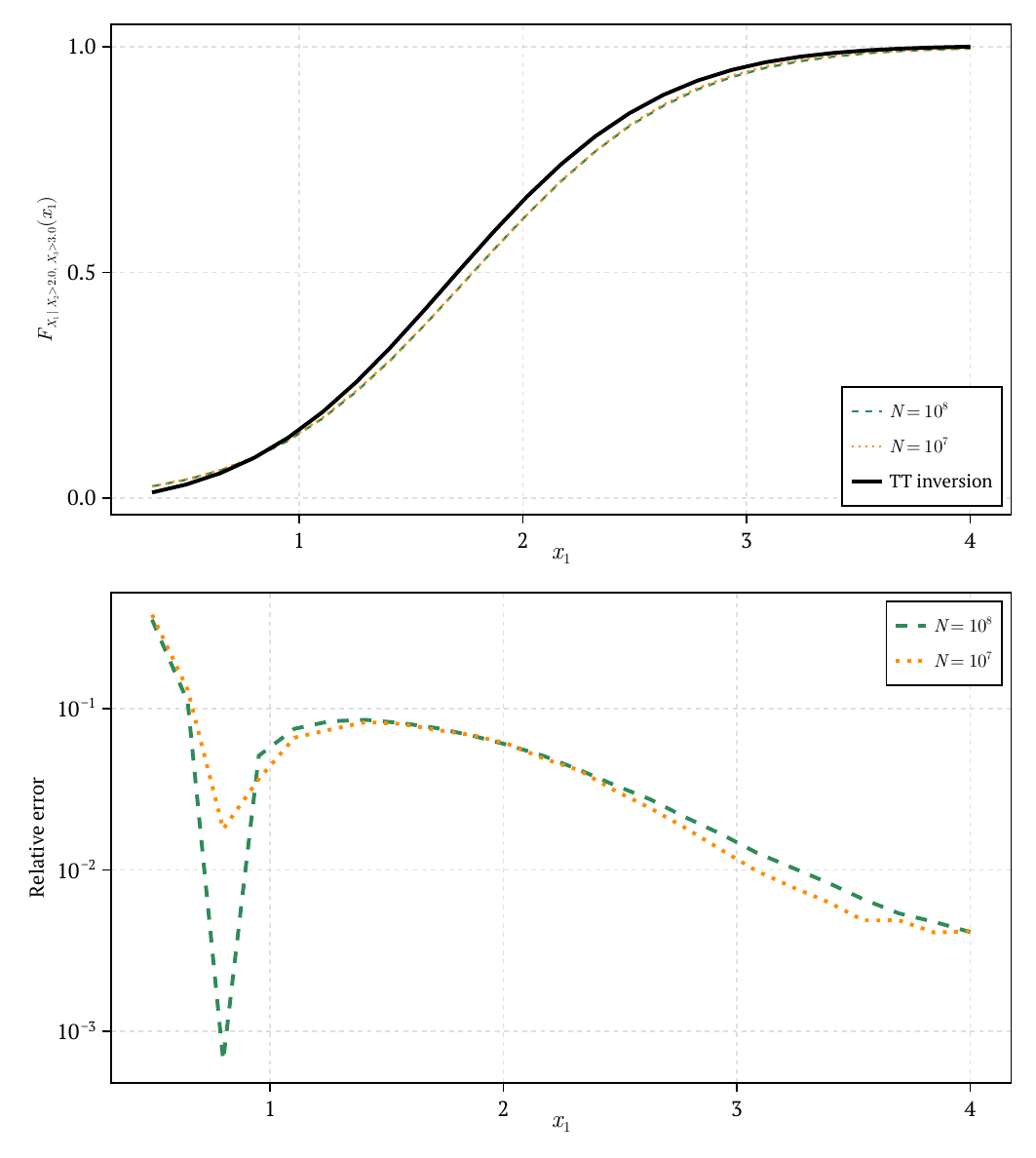}
\caption{Conditional CDF $F_{X_1\mid X_2>2.0,\,X_3>3.0}(x_1)$ computed via TT inversion (black solid), compared to Monte Carlo estimates with $N=10^7$ (orange dotted, min--max band over 10 runs) and $N=10^8$ (green dashed). Inversion parameters $A=15$, $\ell=15$, $m=15$, $n_{\mathrm{E}}=15$, and thresholds $c_2=2.0$, $c_3=3.0$}
\label{fig:conditional}
\end{figure}

Figure~\ref{fig:conditional} compares the TT result against two Monte Carlo estimates obtained by direct simulation from the variance-mean mixture representation. 

\subsection{Mutual information}
As an example of an information-theoretic dependence measure, we compute the pairwise mutual information between components $X_i$ and $X_j$,
\begin{equation}\label{eq:mutual}
I(X_i;\,X_j) = H(X_i) + H(X_j) - H(X_i,X_j),
\end{equation}
where
\begin{equation}
H(X_k) =-\int f_k(t)\log f_k(t) \, \text{d}t
\end{equation}
and
\begin{equation}
H(X_i,X_j) = -\iint f_{ij}(s,t)\log f_{ij}(s,t) \, \text{d}s \,\text{d}t
\end{equation}
denote the differential entropy of the one- and two-dimensional marginals, respectively. In practice, these integrals are evaluated using the same grid-based quadrature as in the normalization and marginal contractions.

\begin{figure}[ht]
\centering
\includegraphics[width=0.5\textwidth]{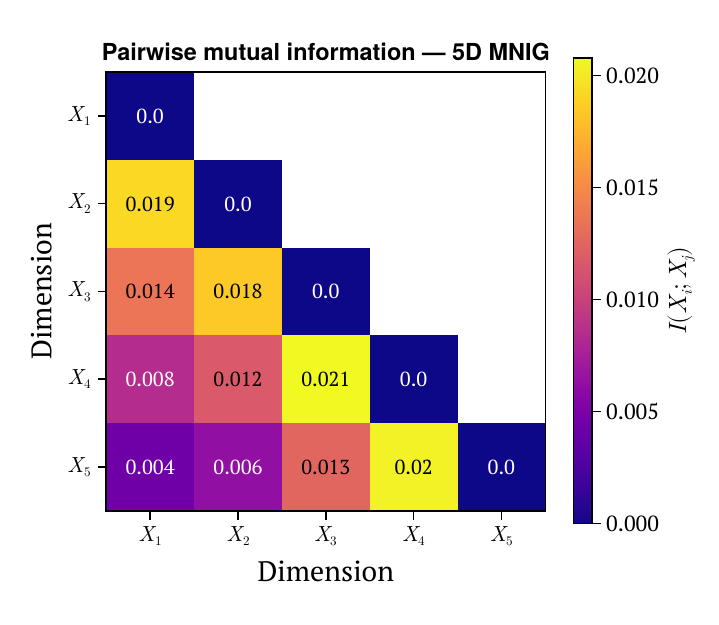}
\caption{Pairwise mutual information $I(X_i \; ;,X_j)$ for the $d=5$ MNIG distribution, computed from the normalized TT representation via quadrature on the one- and two-dimensional marginals. Pairs coupled more strongly through the covariance structure exhibit larger mutual information.}
\label{fig:mutual}
\end{figure}

The resulting mutual information matrix for the MNIG example is shown in Figure~\ref{fig:mutual}. As expected, the strongest dependence occurs for neighbouring components, while pairs that are farther apart in the covariance structure exhibit smaller mutual information.

\section{Discussion}
The main observation of this work is that the multidimensional inverse Laplace transform admits a separable tensor-network structure. In the inversion formula of Choudhury, Lucantoni, and Whitt \cite{MultidimensionalLaplace}, each quadrature node depends only on the local indices $(p_k,j_k,t_k)$ of a single dimension, while all inter-dimensional coupling is contained in the transformed function $\tilde f$, as shown in \eqref{eq:invmultidim}. This separability allows the inversion to be written as a sequence of local tensor contractions acting on the TT representation of $\tilde f$ evaluated at the quadrature nodes \eqref{eq:quad_nodes}. Because the oscillatory factors and Euler weights act independently on each triplet of tensor cores, the contraction does not increase the inter-dimensional bond dimensions; see Figure~\ref{fig:contraction_scheme}.

The numerical examples suggest that several transformed probability distributions admit low-rank representations on the complex quadrature grid. This is particularly useful when the density in the $t$-domain requires expensive evaluations of special functions. In the Wishart and MNIG examples, the maximum bond dimension remains moderate even up to dimension $d=8$, provided the correlations are not too strong. Figures~\ref{fig:MNIG_correlationplot} and \ref{fig:MNIG_scalingplot} show that the bond dimension increases with the correlation parameter $\rho$, reflecting the growing non-separability of $\tilde f$. As $\rho$ approaches one, the covariance matrix becomes nearly uniform, the transformed function becomes less compressible, and the bond dimension grows rapidly. The practical efficiency of the method therefore depends strongly on the compressibility of $\tilde f$. In the examples considered here, the transforms remain sufficiently smooth on the quadrature grid for \texttt{TT-cross} to construct accurate approximations using far fewer function evaluations than the full tensor-product grid.

Our approach differs from methods such as \cite{dolgov2020approximation}, where the probability density itself is approximated directly in TT format. In many applications the density may involve expensive evaluations of special functions or oscillatory integrals, whereas the Laplace transform admits a significantly simpler analytical representation. This is clear for the multivariate normal-inverse Gaussian distribution, where the density contains modified Bessel functions while the Laplace transform is expressed in closed form. By constructing the TT representation from the transform rather than the density, the computational cost is shifted toward evaluations of $\tilde f$ on the contour grid, which can be substantially cheaper whenever the transform retains a simple analytic structure. The resulting TT representation of the recovered density can then be used directly for CDFs, rectangle probabilities, conditional distributions, tail probabilities, marginalisation, entropy estimation, and other distributional queries through tensor contractions.

The method has several limitations. Its computational complexity depends critically on the bond dimension, and strongly correlated multidimensional systems may become too expensive to compress. Similarly, transforms with singularities close to the inversion contour, or with highly oscillatory behaviour, may require finer quadrature grids and careful tuning of the inversion parameters. Although the contraction scheme in Figure~\ref{fig:contraction_scheme} scales linearly in $d$ for fixed $\chi$, the overall cost need not remain linear in regimes where the bond dimension grows rapidly.

Several extensions are possible. One direction is to combine the inversion procedure with TT-based conditional distribution samplers such as \cite{dolgov2020approximation}, enabling posterior sampling from distributions available only through their Laplace transform. Another possibility is to incorporate quantics tensor-train (QTT) representations \cite{qtt_khoromskij}, in which the discretisation itself is encoded hierarchically. Preliminary results suggest that interpolation-based QTT constructions can obtain efficient representations using fewer function evaluations than standard \texttt{TT-cross} \cite{lindsey2023multiscale}. Such an approach could provide finer control over quadrature resolution and enable adaptive or multiscale contour discretisation. In particular, one could represent the quadrature grid in QTT form and then reshape it into a TT for inversion, preserving the separability exploited in this work while improving resolution.

\subsection{Data Availability}

The numerical experiments were performed using \texttt{TensorTrainNumerics.jl} version 1.1.0 \cite{Martin_TensorTrainNumerics_jl_2026}. The implementation used to reproduce the figures and numerical experiments is available at:

\begin{center}
\url{https://github.com/MartinMikkelsen/LaplaceInversion.jl}
\end{center}

\newpage 
\bibliographystyle{plainnat}
\bibliography{references}

\onecolumn 
\appendix

\section{Two-dimensional inversion with continuous variable}\label{app:2D}
The following derivation summarizes \cite{MultidimensionalLaplace}. Let $f(t_1,t_2)$ be a complex-valued function, and let its two-dimensional Laplace transform be
\begin{equation}
    \tilde{f}(s_1,s_2) = \int_0^\infty \int_0^\infty \exp(-(s_1t_1+s_2t_2)) f(t_1,t_2) \, \text{d}t_1 \text{d}t_2,
\end{equation}
assumed to be well-defined. We now exploit the two-dimensional Poisson summation formula
\begin{equation}\label{eq:poisson}
    \sum_{j=-\infty}^{\infty} \sum_{k=-\infty}^\infty F\left( t_1+\frac{2\pi j}{h_1}, t_2+\frac{2\pi k}{h_2} \right) = \sum_{j=-\infty}^{\infty} \sum_{k=-\infty}^\infty \frac{h_1h_2}{4\pi^2} \phi(jh_1, kh_2) \exp(-i(jh_1t_1+kh_2t_2)),
\end{equation}
where $\phi(u_1,u_2)$ is the bivariate Fourier transform given by 
\begin{equation}\label{eq:fourier}
    \phi(u_1,u_2) = \int_{-\infty}^\infty \int_{-\infty}^\infty \exp(-i(t_1u_1+t_2u_2)) F(t_1,t_2) \, \mathrm{d}t_1 \mathrm{d}t_2.
\end{equation}
The main idea is that the left-hand side in equation \eqref{eq:poisson} can be reconstructed by aliasing and the right-hand side can be interpreted as a trapezoidal rule form applied to the Fourier inversion formula given by
\begin{equation}
    F(t_1,t_2) = \frac{1}{4\pi^2} \int_{-\infty}^\infty \int_{-\infty}^{\infty} \exp(-i(t_1u_1+t_2u_2))\phi(u_1,u_2) \, \text{d}u_1 \text{d}u_2
\end{equation}
We introduce exponential damping by making the following substitution for the original function $f$:
\begin{equation}
    F(t_1,t_2) \rightarrow f(t_1,t_2) \exp(-(a_1 t_1+a_2 t_2)), \quad t_1,t_2 \geq 0.
\end{equation}
This allows us to rewrite equation \eqref{eq:poisson} as 
\begin{multline}
    \sum_{j=0}^{\infty} \sum_{k=0}^\infty \exp\!\left(-\left[a_1(1+2jl_1)t_1+a_2(1+2kl_2)t_2\right]\right) f((1+2jl_1)t_1, (1+2kl_2)t_2) \\
    = \frac{1}{4l_1t_1l_2t_2} \sum_{j=-\infty}^{\infty} \sum_{k=-\infty}^{\infty} \exp\left( -i \left(\frac{j\pi}{l_1}+\frac{k\pi}{l_2} \right) \right) \tilde{f}\left(a_1-\frac{ij\pi}{l_1t_1},a_2-\frac{ik\pi}{l_2t_2} \right),
\end{multline}
where we let $a_1=\frac{A_1}{2t_1l_1}$ and $a_2=\frac{A_2}{2t_2l_2}$. This means we arrive at 
\begin{multline}\label{eq:fbar}
    \bar{f}(t_1,t_2)=\frac{\exp(A_1/(2l_1))}{2l_1t_1}\sum_{j=-\infty}^{\infty}\exp\left(-\frac{ij\pi}{l_1}\right) \\
    \left[\frac{\exp(A_2/(2l_2))}{2l_2t_2} \sum_{k=-\infty}^\infty \exp\left( - \frac{ik\pi}{l_2}\right) \left(\tilde{f}\left(\frac{A_1}{2l_1t_1}-\frac{ij\pi}{l_1t_1}, \frac{A_2}{2l_2t_2}-\frac{ik\pi}{l_2t_2} \right) \right)\right],
\end{multline}
where the error term in $f(t_1,t_2)=\bar{f}(t_1,t_2)-\bar e$ has been omitted. We now build an explicit construction of \eqref{eq:fbar} using Euler summation:
\begin{align}\label{eq:EulerSum}
    E(m,n) & = S_n + (-1)^{n+1}\sum_{k=0}^{m-1} (-1)^k 2^{-(k+1)} \Delta^k \alpha_{n+1} \\
    &= \sum_{k=0}^m \binom{m}{k} 2^{-m}S_{n+k}, \quad S_j = \sum_{k=0}^j(-1)^k\alpha_k,
\end{align} 
where $\Delta \alpha_j = \alpha_{j+1}-\alpha_j$ and $\Delta^k$ is the $k$-th application of the forward-difference operator, $\Delta$. Rewriting equation \eqref{eq:fbar} yields

\begin{multline}\label{eq:fbar2}
    \bar{f}(t_1,t_2)=\frac{\exp(A_1/(2l_1))}{2l_1t_1} \sum_{j_1=1}^{l_1} \sum_{j=-\infty}^{\infty} (-1)^j  \exp\left(-\frac{ij_1\pi}{l_1}\right) \\
    \left[\frac{\exp(A_2/(2l_2))}{2l_2t_2}\sum_{k_1=1}^{l_2} \sum_{k=-\infty}^\infty (-1)^k \exp\left( - \frac{ik_1\pi}{l_2}\right) \left(\tilde{f}\left(\frac{A_1}{2l_1t_1}-\frac{ij_1\pi}{l_1t_1}- \frac{ij\pi}{t_1}, \frac{A_2}{2l_2t_2}-\frac{ik_1\pi}{l_2t_2} - \frac{ik\pi}{t_2}\right) \right)\right].
\end{multline}
\subsection{2D validation}
Below we test that each of the examples in \ref{sec:examples} agree in 2D.

\begin{figure}[ht]
\centering
\includegraphics[width=\textwidth]{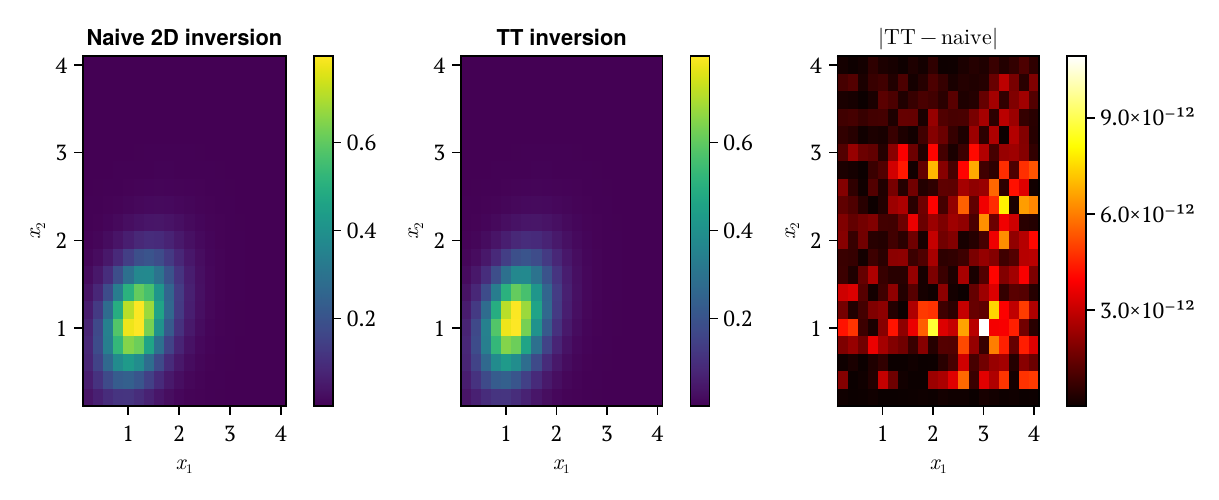}
\caption{The 2D MNIG naïve inversion algorithm}
\label{fig:2DMNIG}
\end{figure}

\begin{figure}[ht]
\centering
\includegraphics[width=\textwidth]{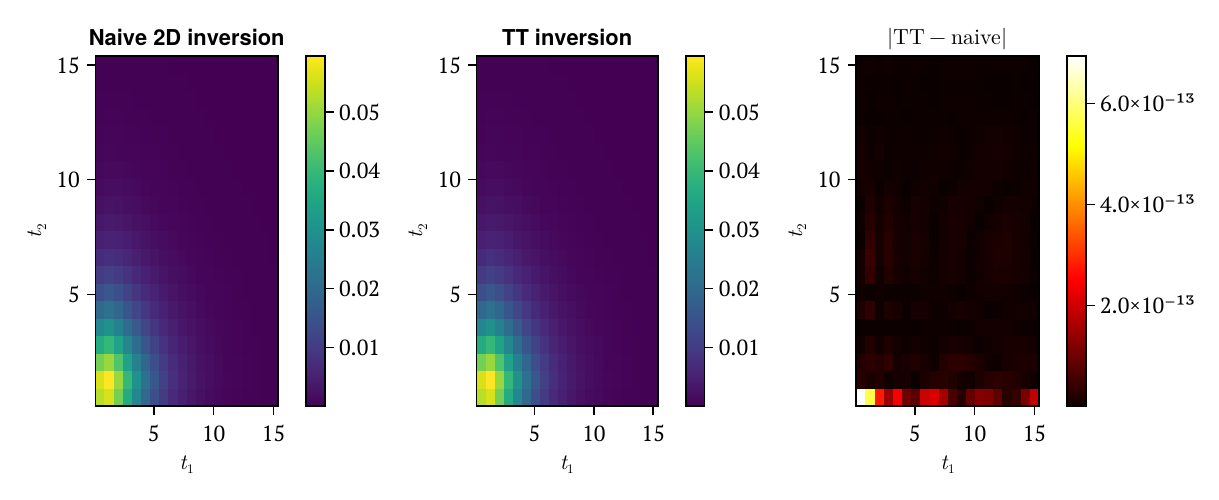}
\caption{The 2D Wishart naïve inversion algorithm}
\label{fig:2DWishart}
\end{figure}

\begin{figure}[ht]
\centering
\includegraphics[width=\textwidth]{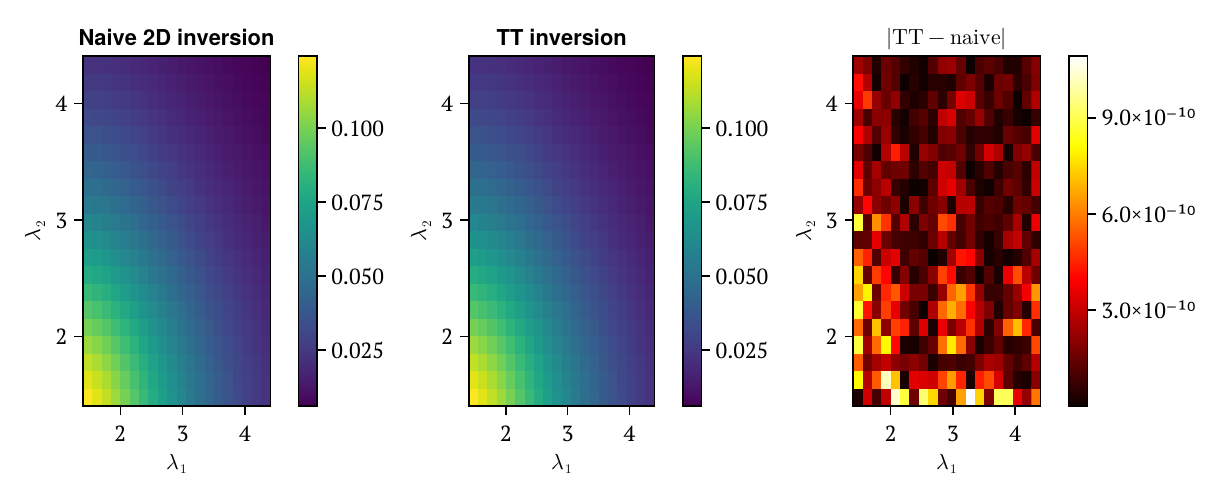}
\caption{The 2D correlated Gamma model compared to the naïve inversion algorithm}
\label{fig:2Dcorrelated}
\end{figure}

\end{document}